\documentclass[report]{svmult}

\usepackage{mathptmx}        
\usepackage{helvet}          
\usepackage{courier}         
\usepackage{type1cm}         

\usepackage{makeidx}         
\usepackage{graphicx}        
\usepackage{multicol}        
\usepackage[bottom]{footmisc}

\usepackage{amsmath,amssymb,bbm}
\usepackage{graphicx}
\usepackage[utf8]{inputenc}
\usepackage{url}
\usepackage{microtype}
\usepackage{import}
\usepackage{color}
\usepackage{epic}
\usepackage{tikz}
\usepackage{etex}
\usepackage{pgfplots}
\usepackage{subfig}
\usepackage[ruled,vlined]{algorithm2e}

\usepackage{graphicx}
\usepackage{bm}

\usepackage{appendix}
\DeclareMathAlphabet{\mathitbf}{OML}{cmm}{b}{it}

\newcommand{\bzero}{\textbf{0}}
\newcommand{\vTheta}{\varTheta}

\newcommand{\Pol}{\Psi}
\newcommand{\pol}{\psi}
\newcommand{\EXP}[1]{\mathbb{E}\left(#1\right)}
\newcommand{\di}{\mathrm{d}}
\newcommand{\btheta}{\theta}

\newcommand{\mysol}{c}
\newcommand{\Psii}{\mysol}

\newcommand{\bsol}{\mathbf{\mysol}}
\newcommand{\Do}{\mathcal{D}}

\newcommand{\var}[1]{{\ensuremath{\mathrm{Var}}\mspace{-2mu}\left[#1\right]}}

\def\sol{c}
\def\bsol{\mathbf{\sol}}

\newcommand{\myug}{ug4\;}

\newcommand{\gPC}{gPC\;}
\newcommand{\QMC}{qMC\;}

\newcommand{\thetab}{{\theta}}

\newcommand{\bx}{\mathbf{x}}

\newcommand{\be}{\mathbf{e}}

\newcommand{\conc}{c} 
\newcommand{\pres}{p} 
\newcommand{\poro}{\phi} 
\newcommand{\perm}{{\mathbf{K}}} 
\newcommand{\dens}{\rho} 
\newcommand{\visc}{\mu} 
\newcommand{\dvel}{{\mathbf{q}}} 
\newcommand{\grav}{{\mathbf{g}}} 
\newcommand{\disp}{{\mathbf{D}}} 
\begin{document}
\graphicspath{{../figures/}} 

\title*{Propagation of Uncertainties in Density-Driven Flow}
\author{Alexander Litvinenko\inst{1} \and Dmitry Logashenko\inst{2} \and Raul Tempone\inst{1,2} \and
Gabriel Wittum\inst{2,3} \and David Keyes\inst{2}}
 \authorrunning{Litvinenko et al.}
\institute{\inst{1}RWTH Aachen, Kackertstr. 9C, Aachen, Germany, Phone: +492418099203, \email{litvinenko, tempone@uq.rwth-aachen.de} \\ \inst{2}KAUST, Thuwal/Jeddah, Saudi Arabia, 
\email{dmitry.logashenko, gabriel.wittum, david.keyes@kaust.edu.sa} \\ \inst{3}G-CSC, Frankfurt University, Kettenhofweg 139, Frankfurt, Germany}
%
%
\maketitle

\abstract*{Accurate modeling of contamination in subsurface flow and water aquifers is crucial for agriculture and environmental protection. Here, we demonstrate an efficient parallel algorithm to quantify the propagation of the uncertainty in the dispersal of pollution in subsurface flow. Specifically, we consider the density-driven flow and estimate how uncertainty from permeability and porosity propagates to the solution. We take an Elder-like problem as a numerical benchmark, and we use random fields to model the limited knowledge on the
porosity and permeability. We construct a low-cost generalized polynomial chaos (\gPC) expansion surrogate model, where the \gPC coefficients are computed by projection on sparse and full tensor grids. We parallelize both the numerical solver for the deterministic problem based on the multigrid method and the quadrature over the parametric space.}

\abstract{Accurate modeling of contamination in subsurface flow and water aquifers is crucial for agriculture and environmental protection. Here, we demonstrate a parallel method to quantify the propagation of the uncertainty in the dispersal of pollution in subsurface flow. Specifically, we consider the density-driven flow and estimate how uncertainty from permeability and porosity propagates to the solution.
We take an Elder-like problem as a numerical benchmark and we use random fields to model the limited knowledge on the
porosity and permeability. We construct a low-cost generalized polynomial chaos expansion (\gPC) surrogate model, where the \gPC coefficients are computed by projection on sparse and full tensor grids. We parallelize both the numerical solver for the deterministic problem based on the multigrid method, and the quadrature over the parametric space.}

\section{Introduction}
\label{litv:sec:Intr}

Accidental contamination of groundwater can be extremely hazardous and thus, accurately predicting the fate of pollutants
in groundwater is essential. Certain pollutants are soluble in water and can leak into groundwater systems, such as seawater into
coastal aquifers or wastewater leaks. Indeed, some pollutants can change the density of a fluid and induce density-driven flows
within the aquifer. This causes faster propagation of the contamination due to convection. Thus, simulation
and analysis of this density-driven flow plays an important role in predicting how pollution can migrate through an aquifer \cite{Kobus2000, Density-driven_Fan97}.

In contrast to transport of pollution by molecular diffusion, convection due to density-driven flow is
an unstable process that can undergo quite complicated patterns of distribution. 
The Elder problem is a simplified but comprehensive model that
describes the intrusion of salt water from a top boundary into an aquifer \cite{Voss_Souza,Elder_Overview17}. The evolution
of the concentration profile is typically referred to as \textit{fingering}. Note that due
to the nonlinear nature of the model, the distribution of the contamination strongly
depends on the model parameters, so that the system may have several stationary solutions
\cite{Johannsen2003}.

Hydrogeological formations typically have a complicated and heterogeneous structure, and geological media may consist of layers of porous media
with different porosities and permeability coefficients \cite{ScheiderKroehnPueschel2012, ReiterLogashenkoVogelWittum2017}. Difficulty in specifying hydrogeological parameters of the media, as
well as measuring the position and configuration of the layers, gives rise to certain errors. Typically, the averaged values of these
quantities are used. However, due to the nonlinearities of
the flow model, averaging of the parameters does not necessarily lead to
the correct mathematical solution. Uncertainties arise from various factors, such as inaccurate measurements of the different parameters and the inability to
measure parameters at each spatial location at any given time. To model the uncertainties, we can use random variables, random fields, and random processes.
However, uncertainties in the input data can spread through the model and therefore also make the solution
uncertain or incorrect. Thus, an accurate estimation of the output uncertainties is crucial, for example, for optimal control and design of the experiment.

Many techniques are used to investigate the propagation of uncertainties associated with porosity and permeability into the
solution, such as classical Monte Carlo sampling.
Although Monte Carlo sampling is dimension-independent, it converges slowly and therefore requires a lot of samples and these simulations may be time-consuming and costly. Alternative, less costly techniques use
collocations, sparse grids and surrogate methods, each with their own advantages and disadvantages. But even advanced
techniques may require hundreds to thousands of simulations and assume a certain smoothness in the quantity of interest. Our simulations
contain up to 0.5-2Mi spatial mesh points and 1000-3000 time-steps, and are therefore run on a massive parallel cluster. Here, we
develop methods that require much fewer simulations, from a few dozen to a few hundred simulations. 

Perturbation methods \cite{CREMER15_Fingers}
are another class of methods that decompose the quantiles of interest (QoIs) with respect to random parameters in a Taylor series.
For small perturbations, the higher order terms can be neglected, thus simplifying the analysis and numerics. These methods assume that random perturbations are small, up to 5$\%$ of the mean, depending on the problem in question. For larger perturbations, these methods are usually not valid.

Here, we use a generalized polynomial chaos expansion (\gPC) technique, where we compute the \gPC coefficients by applying a Clenshaw-Curtis quadrature rule \cite{Barthelmann2000}. We use \gPC to compute QoIs such as the mean, the variance and the exceedance probabilities of the mass fraction (in this case, a saltwater solution). We validate our obtained results using the quasi-Monte Carlo approach. Both methods require the computation of multiple simulations (scenarios) for variable porosity and permeability coefficients. 

To the best of our knowledge, no other reported works have solved the Elder’s problem \cite{Voss_Souza,Elder_Overview17} with uncertain porosity and
permeability parameters using \gPC.
The highly unstable free convective flow has been studied \cite{Xie12_Fingers}, where the number of fingers
and the deepest plume front, the vertical center of solute mass, total solute mass, and solute flux through the source zone was
investigated. The authors concluded that the development of a stochastic framework for modeling free convection is required,
which has significant consequences for model simulation and testing, as well as process prediction. Overviews of the
uncertainties in modeling groundwater solute transport \cite{OverviewUncert93} and modeling soil processes \cite{SoilOverview16} have been performed, as well as
reconnecting stochastic methods with hydrogeological applications \cite{NowakStochMethods18}, which included recommendations for optimization
and risk assessment. Fundamentals of stochastic hydrogeology, an overview of stochastic tools and accounting for uncertainty
have been described in \cite{rubin2003applHydro}.

Recent advances in uncertainty quantification, probabilistic risk assessment, and decision-making under uncertainty in hydrogeologic applications have been reviewed in \cite{TARTAKOVSKYI_Risk13}, where the
author reviewed probabilistic risk assessment methods in hydrogeology under parametric, geologic and model uncertainties.
Density-driven vertical transport of saltwater through the freshwater lens has been modeled in \cite{POST17_Density-driven}.

Various methods can be used to compute the desired statistics, such as direct integration methods (MC, qMC, collocation)
and surrogate-based methods (generalized polynomial chaos approximation, stochastic Galerkin) \cite{Philipp12,matthiesKeese03cmame, babuska2004galerkin}.
Direct integration methods compute statistics by sampling 
unknown input coefficients (they are undefined and therefore uncertain)
and solving the corresponding PDE, while the surrogate-based method computes a low-cost functional (polynomial, exponential, trigonometrical) approximation of QoI.
Examples of the surrogate-based methods include radial basis functions \cite{liu2014,bompard2010,Loeven2007,giunta2004}, sparse polynomials \cite{
chkifa-adapt-stochfem-2015,Sudret_sparsePCE}, polynomial chaos expansion \cite{Xiu, Habib09_PCE, ConradMarzouk13, Dongbin} and low-rank tensor approximation \cite{DolgLitv15, Kim2006,dolgov2014computation,LitvSampling13,Litvin11PAMM}. A nonintrusive stochastic Galerkin was introduced in \cite{Giraldi2014}. The surrogate-based methods have an additional advantage if the gradients are available \cite{liu2017quantification}.  
Sparse grids methods to integrate high-dimensional integrals are considered in \cite{smoljak63, Griebel_Bungartz, Griebel, spiterp, novakRitter97, Novak1996, gerstnerGriebel98-numint, novakRitter99-simple, ConradMarzouk13, CONSTANTINE12}.  A software package is available in \cite{petrasSmolpak}.

The quantification of uncertainties in stochastic PDEs can pose a great challenge due to the possibly large number of random variables involved, and the high cost of each deterministic solution.

For problems with a large number of random variables, the methods based on a regular grid (in stochastic space) are less preferable due to the high computational cost. Instead, methods based on a scattered sampling scheme (MC, qMC) give more freedom in choosing the sample number $N$ and protect against sample failures. The MC quadrature and its variance-reduced variants have dimension-independent convergence $\mathcal{O}(N^{-\frac{1}{2}})$, and qMC has the worst case convergence $\mathcal{O}(\log^M(N)N^{-1})$, where $M$ is the dimension of the stochastic space \cite{matthies2007}. Collocations on sparse or full grids are affected by the dimension of the
integration domain \cite{babuska_collocation, nobile-sg-mc-2015, NobileTemponeWebster08}. 
In this work, we use the Halton rule to generate quasi-Monte Carlo quadrature points \cite{caf1998, joe2008}. A numerical comparison of other qMC sequences has been described in \cite{radovic1996}.

The construction of a \gPC-based surrogate model \cite{Matthies_encicl} is similar to computing the Fourier expansion, where only the Fourier coefficients need to be computed.
Some well-known functions, such as multivariate Legendre, Hermite, Chebyshev and Laguerre functions, are taken as a basis \cite{xiuKarniadakis02a}. These functions should possess some useful properties, e.g., orthogonality. Usually a
small number of quadrature points are used to compute these coefficients. Once the surrogate is constructed, its sampling could
be performed at a lower cost than a sampling of the original stochastic PDE.

To tackle the high numerical complexity, we implement a two-level parallelization.
In the first level, we run all available scenarios in parallel, with each scenario on a separate computing node. On the second level, each scenario is computed in parallel on $2-32$ cores.

This work is structured as follows: In Section \ref{litv:subsec:setup}, we outline the model of the density-driven groundwater flow in porous media and the numerical methods for this type of problem. We describe the stochastic modeling, integration methods and the generalized polynomial chaos expansion technique in Section~\ref{litv:sec:StochModel}.
We present details of the parallelization of the computations in Section~\ref{sec:parallelization}.
Our multiple numerical results, described in Section \ref{litv:sec:numerics}, demonstrate the practical performance of the parallelized solver for the Elder-type problems with uncertainty coefficients in 2D computational domains. We conclude this work with a discussion in Section~\ref{litv:sec:Conclusion}.

\textbf{Our contribution.}  
We applied an established \gPC technique to estimate uncertainties in a density-driven flow. To this end, we solved a time-dependent, nonlinear, second order differential equation with uncertain coefficients. We estimated the propagation of input uncertainties associated with the porosity and permeability into the QoI, which could be a) the mean and the variance of the mass fraction in a sub-domain at each time step, b) the exceedance probability or quantiles in selected points, or c) probability density function in a point. We propose that these statistics could be further used for more efficient Bayesian inference, data assimilation, optimal design of the experiment and optimal control. Our overall goal is to model the dispersion of water pollution and monitor the movement of subsurface water flow. In addition, we used the novel technique of parallelization in spatial and stochastic spaces. We solved each deterministic problem and all stochastic realizations in parallel.

Although we obtained promissing preliminary results using the \gPC technique, many questions remain. For example, a more advanced sparse grid rule may be needed to approximate the \gPC coefficients of the solution.
\section{Problem setup}
\label{litv:subsec:setup}
\subsection{Density-driven groundwater flow problem}
In this section, we summarize the model of the density-driven groundwater flow. For further details, we refer the reader to \cite {Bear-2,Bear1979,Bear-Bachmat,DierschKolditz2002,POST17_Density-driven}.

We consider a domain $\Do \subset \mathbb{R}^d$, $d=2$, filled with two phases: a solid phase and a liquid phase. The solid phase is an immobile solid porous matrix with a mobile liquid in its pores. The matrix is characterized by its porosity $\poro: \Do \to \mathbb{R}$ and its permeability $\perm: \Do \to \mathbb{R}^{d \times d}$. The liquid phase is a solution of salt (sodium chloride) in water. We denote the mass fraction of the brine (a saturated saltwater solution) in the liquid phase by $\conc (t, \mathbf{x}): [0, +\infty) \times \Do \to [0, 1]$. The liquid phase has density $\dens = \dens (\conc)$ and viscosity $\visc = \visc (\conc)$.

The motion of the liquid phase through the solid matrix is characterized by the Darcy velocity $\dvel (t, \mathbf{x}): [0, +\infty) \times \Do \to \mathbb{R}^d$. In these terms, the conservation laws for the liquid phase and the dissolved salt can be written as
\begin{eqnarray}
 \label {e_cont_eq}
 \partial_t (\poro \dens) & + & \nabla \cdot (\dens \dvel) = 0, \\
 \label {e_tran_eq}
 \partial_t (\poro \dens \conc) & + & \nabla \cdot (\dens \conc \dvel - \dens \disp \nabla \conc) = 0,
\end{eqnarray}
where the tensor field $\disp$ represents the molecular diffusion and the mechanical dispersion of the salt. In addition, some momentum equation must be considered for $\dvel$. We assume the Darcy's law for $\dvel$:
\begin{eqnarray} \label {e_Darcy_vel}
 \dvel = - \frac{\perm}{\visc} (\nabla \pres - \dens \grav),
\end{eqnarray}
where $\pres = \pres (t, \mathbf{x}): [0, +\infty) \times \Do \to \mathbb{R}$ denotes the hydrostatic pressure and $\grav = (0, \dots, - g)^T \in \mathbb{R}^d$ is the gravity vector.

For simplicity, we assume an isotropic porous medium characterized by a scalar permeability
\begin{eqnarray} \label {e_scalar_perm}
 \perm = K \mathbf{I},\;\; \text{where}\;\; K = K (\mathbf{x}) \in \mathbb{R},\;\;\text{and}\;\;\mathbf{I}\in\mathbb{R}^{d\times d}\;\text{is the identity matrix}.
\end{eqnarray}
We use the linear dependence for the density:
\begin {eqnarray} \label {e_lin_density}
 \dens (\conc) = \dens_0 + (\dens_1 - \dens_0) \conc,
\end {eqnarray}
where $\dens_0$ and $\dens_1$ denote the densities of pure water and the brine, respectively. Thus, $\conc \in [0, 1]$, where $\conc = 0$ corresponds to the pure water and $\conc = 1$ to the saturated solution. Furthermore, we neglect the dispersion and assume that
\begin {eqnarray} \label {e_mol_diff}
 \disp = \poro D_m \mathbf{I},
\end {eqnarray}
where $D_m$ is the coefficient of the molecular diffusion. The viscosity is considered to be constant. Fields $\phi(\mathbf{x})$ and $K(\mathbf{x})$ will depend on the stochastic variables. Values of the other parameters used in this work are presented in Table \ref{tab:ElderParam}.
\begin{table}[b]
\begin{center}
 \caption{Parameters of the density-driven flow problem}
 \label{tab:ElderParam}
 \begin{tabular}{|l|l|l|} \hline
  Parameters & Values and Units & Description \\ \hline
  $\EXP{\phi}$ & 0.1 [-] & mean value of the porosity \\ \hline
  $D_m$ & $0.565\cdot 10^{-6}$ [$\mathrm{m}^2 \cdot \mathrm{s}^{-1}$] & molecular diffusion \\ \hline
  $\EXP{\perm}$ & $4.845\cdot 10^{-13}$ [$\mathrm{m}^2$] & mean value of the permeability \\ \hline
  $\grav$ & $9.81\be_z$ [$\mathrm{m} \cdot \mathrm{s}^{-2}$] & gravity \\ \hline
  $\rho_0$ & $1000$ [$\mathrm{kg} \cdot \mathrm{m}^{-3}$] & density of pure water \\ \hline
  $\rho_1$ & $1200$ [$\mathrm{kg} \cdot \mathrm{m}^{-3}$] & density of brine \\ \hline
  $\mu$ & $10^{-3}$ [$\mathrm{kg} \cdot \mathrm{m}^{-1} \cdot s^{-1}$] & viscosity \\ \hline
 \end{tabular}
\end{center}
\end{table}

The problem (\ref {e_cont_eq}--\ref {e_mol_diff}) must be closed by the specification for the boundary conditions for equations (\ref {e_cont_eq}) and (\ref{e_tran_eq}), i.e. for $\conc$ and $\pres$, as well as the initial conditions for $\conc$. In our tests, we consider variations of the Elder problem \cite{Voss_Souza,Elder_Overview17}: a heavily concentrated solution intrudes due to gravitation through the upper boundary into an aquifer filled initially with pure water. During this process, a complicated flow arises leading to a specific distribution of the intruding solution, in a process usually described as ``fingering'' \cite{Johannsen2003}. The time evolution of the solution in (\ref{e_cont_eq}--\ref{e_Darcy_vel}) is determined by the initial conditions for $\conc$. However, the system (\ref {e_cont_eq}--\ref {e_Darcy_vel}) is nonlinear, thus the model can have several stationary states \cite{Johannsen2003}. In addition, for certain parameter settings, the flow could be very sensitive to small variations of the initial data. Small fluctuations may lead to a principally different asymptotic distribution of the mass fraction. Spatial variations of the porosity and permeability fields make the situation more complicated. The same initial and boundary data may correspond to principally different asymptotic behaviors of the solution for particular distributions of the parameters of the porous medium. The influence of this factor can be estimated by the application of the stochastic models.

\begin{figure}[htbp!]
\centerline{\includegraphics[width=0.7\textwidth]{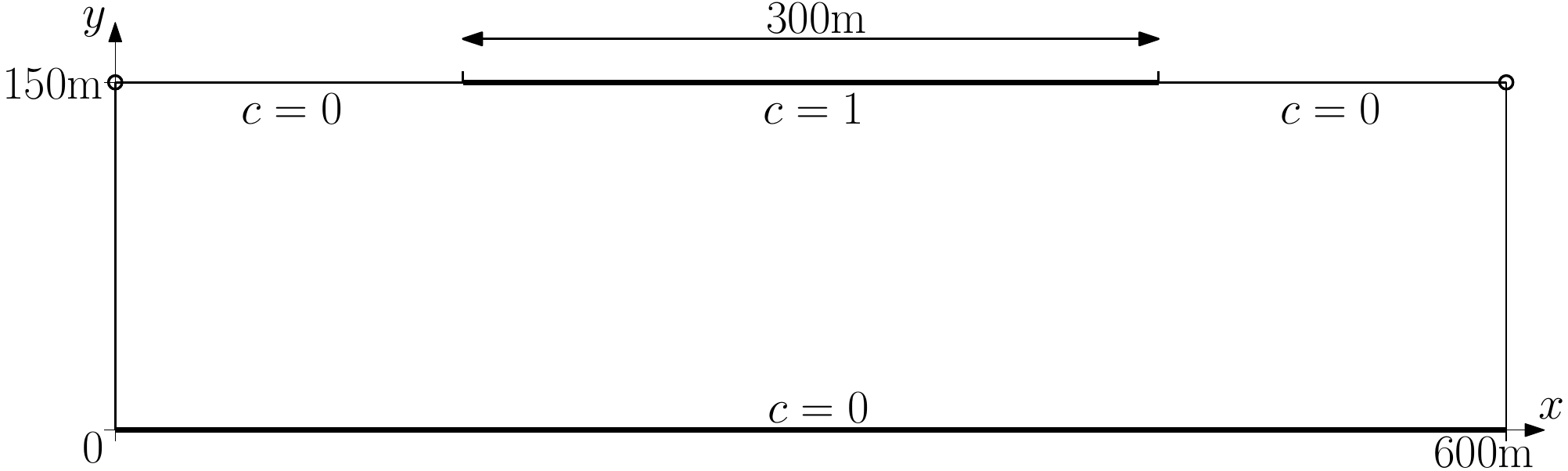}}
\caption{2D reservoir geometry $\Do = (0,600)\times(0,150)$ $\mathrm{m}^2$.}
\label{fig:Elder2}
\end{figure}

In this work we consider a 2D reservoir geometry $\Do = (0,600)\times(0,150)$ $\mathrm{m}^2$ (cf.\ Fig.~\ref{fig:Elder2}). We set the Dirichlet boundary conditions for $\mysol$ on the top and bottom boundaries according to this scheme. On the left and right sides, we impose the no-flux boundary conditions for equation (\ref {e_tran_eq}). For equation (\ref{e_cont_eq}), we impose the no-flux boundary conditions on the whole $\partial \Do$. However, to fix $p$, we set it to $0$ at the two upper corners of the domain (denoted by the small circles on the scheme).
\subsection{Modeling of porosity, permeability and mass fraction}
\label{litv:subsec:PorosityVar}
We model unknown input parameters - the porosity ($\poro$) and permeability ($\perm$) of the solid phase by random fields.
The mass fraction $\mysol$ is the unknown output QoI, and
is a function of $\poro$ and $\perm$.
We introduce a randomness in $\poro$ and assume $\perm$ to be isotropic and dependent on $\poro$:
\begin {eqnarray} \label {e_perm_of_poro}
 \perm = K \mathbf{I}, \qquad \text{and}\quad K = K (\poro). 
\end {eqnarray}
The distribution of $\poro(\bx,\thetab)$, $\bx\in \Do$, is determined by a set of stochastic parameters $\btheta=(\theta_1,...,\theta_M,...)$. Each component $\theta_i$ is a random variable.
The dependence (\ref{e_perm_of_poro}) is specific for every material and there is no a general law. In our experiments, we use a Kozeny-Carman-like dependence \cite{Costa_2006}:
\begin{eqnarray} 
\label{e_perm_Kozeny_Carman}
 K (\poro) = \kappa_{KC} \cdot \dfrac {\poro^3} {1 - \poro^2},
\end {eqnarray}
where the scaling factor $\kappa_{KC}$ is chosen to satisfy the equality $K(\EXP{\poro}) \mathbf{I} = \EXP{\mathbf{K}}$. Further parameters of the standard Elder problem are listed in Table \ref{tab:ElderParam}).\\
%
\textbf{Quantities of interest.} The random variability in the porosity may result in very different mass fraction distributions. The number of fingers, the velocity and time of mass fraction propagation in the media and the form and intensity may vary. We compute the mean value and the variance of the mass fraction. The regions with high variance indicate regions with a high variability/uncertainty. Such regions may require additional attention from specialists (a finer mesh or additional observations). For $n_p$ pre-selected points in $\Do$, we evaluate the probability density function and exceedance probabilities. Such probabilities may help to monitor the quality of drinking water and to estimate the risks of contamination.
\subsection{Numerical solution of the flow model}
\label{sec:NumFlow}
For the estimation of the uncertainty, we computed realizations of solutions (\ref{e_cont_eq}--\ref{e_tran_eq})
for sets of the stochastic parameters $\btheta$. We used plugin d${}^3$f of the simulation framework \myug,
\cite{ug4_ref1_2013,ug4_ref2_2013}. This framework presents a flexible tool for the numerical solution of non-stationary and nonlinear systems of PDEs, and is parallelized and highly scalable.
For the solution, Equations~\ref{e_cont_eq}--\ref{e_tran_eq} are discretized by a vertex-centered finite-volume method on
unstructured grids in the geometric space, as presented in \cite{Frolkovic-MaxPrinciple}. In particular, the consistent
velocity approach is used for the approximation of the Darcy velocity \eqref{e_Darcy_vel},
\cite{Frolkovic-ConsVel,Frolkovic-Knaber-ConsVel}. The implicit Euler scheme is used for the time discretization.
The implicit time-stepping scheme is chosen for its unconditional stability, as the velocity depends on the
unknowns of the system. This is especially important for variable coefficients as it is difficult to predict the
range of the variation of the velocity in the realizations.

In this discretization, we obtain a large sparse system of nonlinear algebraic equations in every time step, and we solve it using the Newton method with a line search method. The sparse linear system, which appears in the nonlinear
iterations, is solved by the BiCGStab method (\cite {Templates}) with the multigrid preconditioning (V-cycle,
\cite{Hackbusch85}), which proved to be very efficient in this case. In the multigrid cycle, we use the
ILU${}_\beta$-smoothers \cite{Hackbusch_Iter_Sol}.

\section{Stochastic Modeling and Methods}
\label{litv:sec:StochModel}
In \cite{CREMER15_Fingers}, the authors compare six perturbation methods under saturated and unsaturated flow conditions to generate fingers in laboratory sand tank experiments. They showed that it is possible to reproduce dense plume fingering for saturated flows.
Additionally, the authors studied various perturbation mechanisms to induce instabilities in the numerical solution of variable-density systems.
In this section, we review well-known basic theories and generalized polynomial chaos surrogates.
\subsection{Sampling Methods}
\label{litv:sec:methods}


Let $\btheta=(\theta_1,\ldots,\btheta_M)$ be a vector of independent random variables.
We begin by defining a probability space ${S}:=(\vTheta, \mathcal{B}, \mathbb{P})$, where
$\vTheta$ denotes a sample space, $\mathcal{B}$ a $\sigma$-algebra on $\vTheta$, and $\mathbb{P}$ a probability measure on $(\vTheta, \mathcal{B})$. 
As stated in \cite{wiener38}, any RV $\theta:\vTheta \rightarrow \mathbb{R}^M$ with finite variance can be represented in terms of the polynomial chaos expansion, e.g. as a multi-variate Hermite polynomial of Gaussian RVs. Alternatively to Hermite polynomials, other orthogonal polynomials (Legendre, Chebyshev, Laguerre) \cite{Askey85, xiuKarniadakis02a, Xiu}, or splines, for instance, can be used.
We assume $\{\Pol_{\beta} \}_{\beta \in \mathcal{J}}$ is a basis of $S=L_2(\vTheta, \mathbb{P})$. The cardinality of the multi-index set $\mathcal{J}$ is infinite, therefore, we truncate it and obtain a finite set $\mathcal{J}_{M,p}$ and a subspace $S_{M,p}$. The solution $\mysol(t,\bx,\btheta)$ lies in the tensor product space $L_2(D)\otimes L_2(\vTheta)\otimes L_2([0,T])$, where $[0,T]$ is a time interval, where the initial problem is solved. After a discretization, we obtain $\bsol \in \mathbb{R}^{n}\otimes {S}_{M,P}\otimes \mathbb{R}^{n_t}$, where $n$ is the number of basis functions in the physical space, $\vert \mathcal{J}_{M,p} \vert$ dimension of the stochastic space, and $n_t$ dimension of the temporal space.

The mean value of $\mysol(t,\bx,\btheta)$ can be computed as an integral over the multidimensional domain $\vTheta$:
\begin{equation} \label{eq:pdestat}
  \overline{\mysol}(t,\bx)=\EXP{\mysol(t,\bx,\btheta)} = \int_\vTheta \mysol(t,\bx,\btheta)\,\mathbb{P}(\di \btheta),
\end{equation}
where $\mathbb{P}(\di \btheta)$ is a probability measure.
This integral can be computed numerically with some sparse or dense quadrature rules with points $\theta_i$ and weights $w_i$:
\begin{equation} \label{eq:numer_int}
  \EXP{{\mysol}}=\overline{\mysol}(t,\bx)\approx \sum_{i=1}^{N_q} w_i 
  \mysol(t,\bx, \btheta_i)=\sum_{i=1}^{N_q} w_i \mysol(t,\bx,\btheta_i),
\end{equation}
where $\mysol(t,\bx,\btheta_i)$ is the approximate solution produced by the deterministic solver for the realization
$\btheta_i \in \mathbb{R}^M$.
The formula for the sampled variance is:
\begin{equation} \label{eq:numer_var}
  \var{\mysol}(t,\bx)\approx \sum_{i=1}^{N_q} w_i \left ( \overline{\mysol}(t,\bx) - \mysol(t,\bx,\btheta_i) \right)^2.
\end{equation}

%
%
Other well-known methods to compute Equations~\ref{eq:numer_int} and \ref{eq:numer_var} are quasi-Monte Carlo, Monte Carlo and multi-level Monte Carlo methods.
\subsection{Generalized Polynomial Chaos based surrogate model}
\label{litv:sec:gPCE}
As outlined in \cite{Dongbin, xiuKarniadakis02a, wiener38}, the random variable (field) $\Psii(t,\bx,\thetab)$ can be represented as a series of multivariate Legendre polynomials in uncorrelated and independent uniformly distributed random variables
 $\thetab=(\theta_1,...,\theta_m,...)$:
\begin{equation}
\label{eq:PCEdecom}
\Psii(t,\bx,\thetab)=\sum_{\beta \in \mathcal{J}}\Psii_{\beta}(t,\bx)\Pol_{\beta}(\thetab),
\end{equation}
where $\{\Pol_{\beta}\}$ is a multivariate Legendre basis, $\beta = (\beta_1 ,..., \beta_j ,...)$ is a multiindex and $\mathcal{J}$ is a multiindex set. 

This expansion is called the \textit{generalized polynomial chaos expansion} (\gPC) (see more in the Appendix or in  
\cite{wiener38, xiuKarniadakis02a, Xiu, ernst_mugler_starkloff_ullmann_2012}). 
For numerical purposes, we truncate this infinite \gPC series and obtain a finite series approximation, i.e., we keep only $M$ random variables, $\thetab=(\theta_1,...,\theta_M)$, and limit the maximal order of the multi-variate polynomial $\Pol_{\beta}(\thetab)$ by $p$. We denote the new multiindex subset as $\mathcal{J}_{M,p} \subset \mathcal{J}$. Then the unknown random field $\Psii(t,\bx,\thetab)$ can be approximated as follows:
\begin{equation}
\label{eq:PCEdecom2}
\Psii(t,\bx,\thetab)\approx \widehat{\Psii}(t,\bx,\vTheta)=\sum_{\beta \in \mathcal{J}_{M,p}}\Psii_{\beta}(t,\bx)\Pol_{\beta}(\thetab).
\end{equation}
Here $\Pol_{\beta}(\theta)$ are multi-variate Legendre polynomials defined as follows
\begin{equation*}
  \label{eq:B:def-mult-Legendre}
  \Pol_\beta(\thetab) := \prod_{j=1}^{M} \pol_{\beta_j}(\theta_j);
  \quad \forall \thetab\in\mathbb{R}^M,\, \beta \in\mathcal{J}_{M,p},\;\sum\beta_j \leq p, 
\end{equation*}
$\pol_{\beta_j}(\cdot)$ are Legendre monomials, defined in Eq.~\ref{eq:LegendreMonom}. The scalar product is 
\begin{equation}
\label{eq:spro}
\langle \Pol_{\alpha},\Pol_{\beta} \rangle_{L_2(\vTheta)} = \EXP{\Pol_{\alpha}(\thetab)\Pol_{\beta}(\thetab)}=\int_{\Theta}\Pol_{\alpha}(\thetab)\Pol_{\beta}(\thetab)d\mathbb{P}{\thetab}=Q_{\alpha}\delta_{\alpha\beta},
\end{equation}
where $Q_{\alpha}=\EXP{\Pol_{\alpha}^2}=q_{\alpha_1}\cdot...\cdot q_{\alpha_d}$, $q_{\alpha_j}=\langle \pol_{\alpha_j},\pol_{\alpha_j} \rangle=1/(2\alpha_j +1)$, are the normalization constants and $\delta_{\alpha\beta}=\delta_{\alpha_1\beta_1}\cdot ... \cdot \delta_{\alpha_M \beta_M}$ is the $M$-dimensional Kronecker delta function.

Decomposition \eqref{eq:PCEdecom} can be understood as a response surface for $\Psii(t,x,\thetab)$. As soon as the response surface is built, the value $\Psii(t,x,\btheta)$ can be evaluated for any $\btheta$ almost for free (only by evaluating the polynomial (\ref{eq:PCEdecom})). It can be very practical if $10^6$ samples of $\Psii(t,x,\btheta)$ are needed, for example.
Since Legendre polynomials are $L_2$-orthogonal, the coefficients $\Psii_{\beta}(t,x)$ can be computed by projection: 
\begin{equation}
\label{eq:PCEcoef}
\Psii_{\beta}(t,x)=\frac{1}{\langle \Pol_{\beta},\Pol_{\beta} \rangle}\int_{\Theta}\widehat{\Psii}(t,x,\btheta)\Pol_{\beta}(\btheta)\,\mathbb{P}(\di \btheta)\approx \frac{1}{\langle \Pol_{\beta},\Pol_{\beta} \rangle}\sum_{i=1}^{N_q}\Pol_{\beta}(t,\thetab_i)\widehat{\Psii}(t,x,\thetab_i)w_i,
\end{equation}
where $w_i$ are weights and $\thetab_i$ quadrature points, defined, for instance, by a Gauss-Legendre integration rule.
Once all \gPC coefficients are computed we substitute them into Eq.~\ref{eq:PCEdecom}. For estimates of truncation and aliasing errors, see \cite{ConradMarzouk13}.
There are different strategies to truncate the multi-index set $\mathcal{J}$ to $\mathcal{J}_{M,p}$: a) $\prod_{i=1}^M p_i \leq p$ or b) $\sum_{i=1}^M p_i \leq p$ or c) $p_i\leq p$, $\forall\; i=1..M$. Here $p$ is the multi-variate polynomial degree.
In the multi-index set $\mathcal{J}_{M,p}$ we keep only $M$ RVs, i.e., $\thetab=(\theta_1,\ldots,\theta_M)$ and $p_i$ is the monomial degree w.r.t. the variable $\theta_i$.

\begin{remark}
If only the variance and the mean values are of interest, then we recommend to use the collocation method with some sparse or full-tensor grid. If some additional detail is required, such as sensitivity analysis, parameter identification, or computing cdf or pdf, then we recommend to compute the gPC expansion.
\end{remark}

Quadrature grids for computing multi-dimensional integrals can be constructed from products of 1D integration rules. To keep computational costs to a minimum, a sparse grid technique can be used for certain class of smooth functions \cite{Griebel_Bungartz}. To compute Eq.~\ref{eq:PCEcoef}, we apply the well-known Clenshaw Curtis quadrature rules \cite{CC60}. 
%

%
%
%
Since $\Pol_{\bzero}(\btheta)=1$, it becomes apparent that the mean value is the first \gPC  coefficient 
\begin{equation}
\label{eq:PCEcoef2}
\overline{\sol}:=\EXP{\Psii(t,x,\btheta)}=\Psii_{\bzero}(t,x)=
\frac{1}{\langle \Pol_{\bzero},\Pol_{\bzero} \rangle}\int_{\Theta}{\Psii}(t,x,\btheta)\Pol_{\bzero}(\btheta)\,\mathbb{P}(\di \btheta)
\approx \frac{1}{1}\sum_{i=1}^{N_q}{\Psii}(t,x,\thetab_i)w_i,\quad \bzero=(0,\ldots,0).
\end{equation}

The variance is the sum of squared \gPC coefficients \cite{xiuKarniadakis02a, GhanemBook17,Knio10}:
\begin{align}
\label{eq:var}
\var{\Psii(t,x,\btheta)}&=\EXP{\left(\Psii(t,x,\btheta)-\overline{\Psii}\right) \otimes \left(\Psii(t,x,\btheta)-\overline{\Psii}\right)}=
\sum_{\beta \in \mathcal{J}, \beta>0}\sum_{\gamma \in \mathcal{J}, \gamma>0} \EXP{\Pol_{\beta}(\thetab)\Pol_{\gamma}(\thetab)}\Psii_{\beta}(t,x)\otimes \Psii_{\gamma}(t,x)\\
&=\sum_{\gamma>0}\Vert \Pol_{\gamma} \Vert^2 \Psii_{\gamma}(t,x)\otimes \Psii_{\gamma}(t,x)=
\sum_{\gamma>0}\prod_{k=1}^M\frac{1}{2\gamma_k+1} \Psii_{\gamma}(t,x)\otimes \Psii_{\gamma}(t,x),
\end{align}
where $\otimes$ is a Kronecker product.
\subsection{Computing probability density functions}
After we calculated the \gPC surrogate $\widehat{\sol}(t,\bx,\thetab)$ in Eq.~\ref{eq:PCEdecom2}, we can compute the probability density functions (see Fig.~\ref{fig:pdfs}), exceedance probabilities, and quantilies (see Tab.~\ref{tab:quantiles}) in the selected points.

One should evaluate the multi-variate polynomial $\widehat{\sol}(t,\bx,\thetab)$ in a sufficiently large number of points $\thetab$. This evaluations are of low cost. Later, in Sect.~\ref{sec:pdf} we compute pdf in the given point.

Approximation of the probability density function in a point $(t^{*},\bx^{*})$ is computed by sampling the multivariate polynomial on the right-hand side in Eq.~\ref{eq:PCEpoint}. 
\begin{equation}
\label{eq:PCEpoint}
\widehat{\mysol}(\thetab):=\widehat{\mysol}(t^{*},\bx^{*},\thetab)=\sum_{\beta \in \mathcal{J}_{M,p}}\mysol_{\beta}(t^{*},\bx^{*})\Pol_{\beta}(\thetab).
\end{equation}
The random variable $\widehat{\mysol}(\thetab)$ can be sampled, e.g., $N_s=10^6$ times (no additional extensive simulations are required). The obtained sample can be used to evaluate the required statistics, and the probability density function. The exceedance probability for some critical value $c^{*}$ can be estimated as follows:
\begin{equation}
    P(\mysol>\mysol^{*})\approx \frac{\#\{\widehat{\mysol}(\thetab_i):\; \widehat{\mysol}(\thetab_i)>\mysol^{*}, \; i=1,\ldots,N_s\}}{N_s}.
\end{equation}
Having large enough sample set it is also possible estimate quantiles (see Sect.~\ref{sec:quant}).
We remind that for a random variable $\xi$ and for a fixed value $\alpha\in (0,1)$, the $\alpha$-quantile is the number $\xi_{\alpha}\in \mathbb{R}$ such that $\mathbb{P}(X\leq \xi_{\alpha})\geq \alpha$.
\subsection{Truncation and approximation errors}
\label{sec:trunc}
Analysis of the truncation and approximation errors of \gPC has been studied in, e.g. \cite{CONSTANTINE12, Sinsbeck_2015,ConradMarzouk13}.
From Eq.~\ref{eq:PCEdecom2} it follows that the truncation error is
\begin{equation}
\label{tr_err}
E_t=\Vert \Psii(t,\bx,\thetab)- \widehat{\Psii}(t,\bx,\thetab)\Vert =
\Big \| \sum_{\beta \in \mathcal{J}_{c}}\Psii_{\beta}(t,\bx)\Pol_{\beta}(\thetab) \Big \|,  \quad \mathcal{J}_{M,p} \cup \mathcal{J}_{c}=\mathcal{J}.  
\end{equation}
Additionally, the coefficients $\Psii_{\beta}(t,\bx)$ are also approximated by $\widehat{\Psii}_{\beta}(t,\bx)$. This gives the approximation error
\begin{equation}
\label{eq:appr_err}
E_a=\Big \Vert
\sum_{\beta \in \mathcal{J}_{M,p}}\Psii_{\beta}(t,\bx)\Pol_{\beta}(\thetab)
-
\sum_{\beta \in \mathcal{J}_{M,p}}\widehat{\Psii}_{\beta}(t,\bx)\Pol_{\beta}(\thetab) \Big \Vert
= \Big \Vert \sum_{\beta \in \mathcal{J}_{M,p}}(\widehat{\Psii}_{\beta}(t,\bx)-\Psii_{\beta}(t,\bx))\Pol_{\beta}(\thetab) \Big \Vert
\end{equation}
Summarizing the effect of both errors, we obtain
\begin{equation}
\label{eq:sum_errs}
E_t+E_a=\underbrace{\Big \Vert \sum_{\beta \in \mathcal{J}_c}{\Psii}_{\beta}(t,\bx)\Pol_{\beta}(\thetab) \Big \Vert}_\text{truncation error} +\underbrace{\Big \Vert \sum_{\beta \in \mathcal{J}_{M,p}}(\Psii_{\beta}(t,\bx)-\widehat{\Psii}_{\beta}(t,\bx))\Pol_{\beta}(\thetab) \Big \Vert}_\text{approximation error}.
\end{equation}
\begin{remark}
Spatial resolution plays an important role in the simulations. Thus, the distribution of the parameters of the porous medium should be represented correctly. Furthermore, smoothing the flow on a coarse grid due to numerical diffusion can lead to a loss of some phenomena and therefore reduce the accuracy of the uncertainty quantification. Thus we must cover $\Do$ with sufficiently fine grids and use small time steps. To achieve reasonable computation times, this also assumes parallelization of the evaluation of every realization.
\end{remark}

%
\section{Parallelization}
\label{sec:parallelization}

The computation of the stochastic properties can, therefore, be divided into two phases: computation of the realizations of the flow model at the points of the stochastic grid, and the numerical quadrature over the stochastic domain. For the computations presented in this work, we parallelized both phases. We performed computations on a Cray XC40 parallel cluster Shaheen II provided by the King Abdullah University of Science and Technology (KAUST) in Saudi Arabia. The cluster has 6174 nodes with 2 Intel Haswell 2.3 GHz CPUs per node, and 16 processor cores per CPU. Every node has 128 GB RAM. The nodes communicate via the Cray Aries interconnect network with Dragonfly topology. The system is equipped with the Sonexion Lustre 2000 storage appliance.

In the computations, we computed all the realizations concurrently, and the results for the essential time steps were stored in the parallel file system.

Parallelization of the framework \myug is based on the distribution of the spatial grid between the processes and is implemented using MPI \cite{ug4_ref1_2013}. As the computation of the realizations involves the nonlinear iterations and the solution of the large sparse linear systems (see Section \ref{sec:NumFlow}), quite intensive communication is required. We performed the computation of every realization on the 32 cores of one node. We computed different realizations on different nodes. Shaheen II used the SLURM system to manage parallel jobs, and the job arrays were used for the concurrent computations of the realizations. This was the most time-consuming phase, requiring hours of computation time.

Following computation of the realization, we computed the quadrature on 32 cores of one node. The saved data were from the file system. As the quadrature rules were local in the geometric space, they do not require much communication. 
The computation time for this quadrature required minutes, i.e., much shorter than the computation of the realizations desribed above.
\section{Numerical Experiments}
\label{litv:sec:numerics}
We performed several numerical experiments with $M=3,4,5$ random variables. The coefficients of \gPC were computed with full and sparse Gauss-Legendre/Clenshaw-Curtis tensor grids \cite{CONSTANTINE12}. 
The maximal polynomial order $p$ lay between 2 and 5, depending on the porosity function.
%
%
%

\subsection{Experiment 1. Elder's problem with 3 RVs}
\begin{figure}[htbp!]
    \centering
    \includegraphics[width=0.49\textwidth]{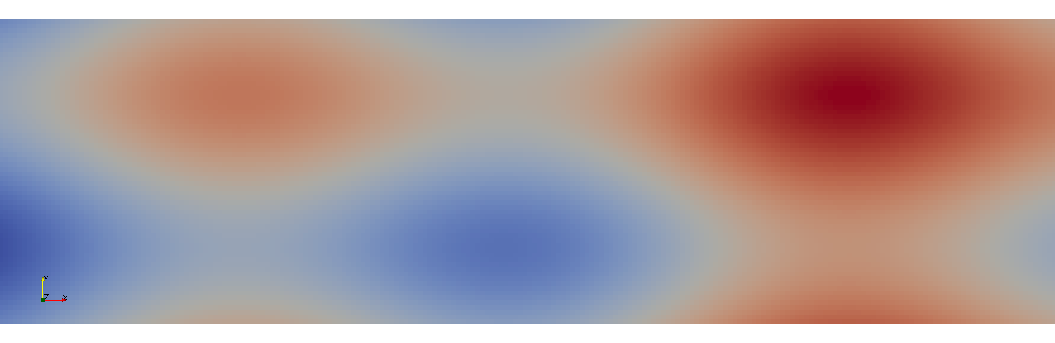}
    \caption{A realization of the porosity field (\ref{eq:poro3}).}
    \label{fig:Ex1-porosity}
\end{figure}

In this experiment, we solve the problem from Section \ref{litv:subsec:setup} (cf.\ Fig.\ \ref{fig:Elder2}) with the following smooth porosity field
\begin{equation}
\label{eq:poro3}
\poro(\bx,\thetab)= 0.1+0.005({\theta_1}\cos(\pi x/600) + \theta_2\sin(2\pi y/150)+ \theta_3 \cos(4\pi x/600)), \quad \bx = (x,y) \in [0,600]\times[0,150],
\end{equation}
parametrized with 3 stochastic parameters $\theta_1,\theta_2,\theta_3 \sim U[-1,1]$. One realization of $\poro(\bx,\thetab)$ is shown in Fig.\ \ref{fig:Ex1-porosity}.

\subsubsection{Dependence of the solution on the spatial and temporal grid.}
\begin{figure}[htbp!]
\centering
\subfloat[]{\includegraphics[width=0.49\textwidth]{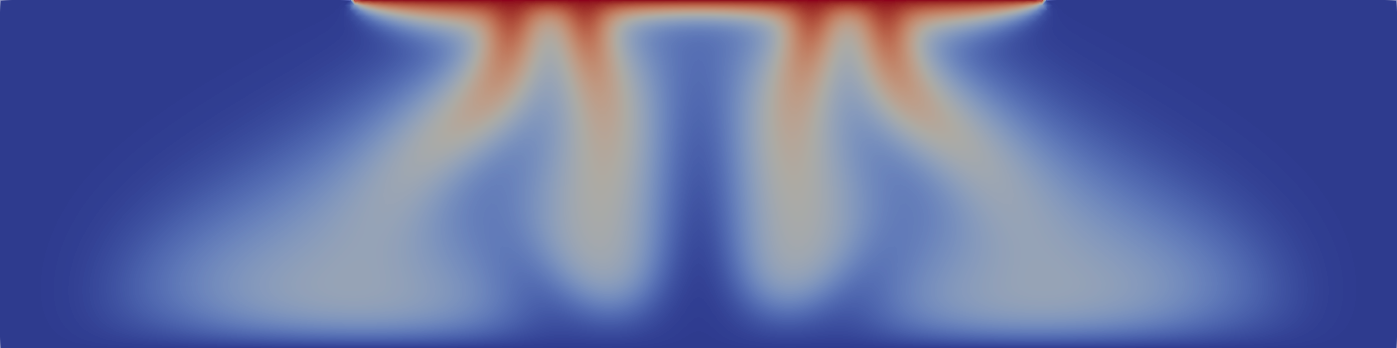}}\,
\subfloat[]{\includegraphics[width=0.49\textwidth]{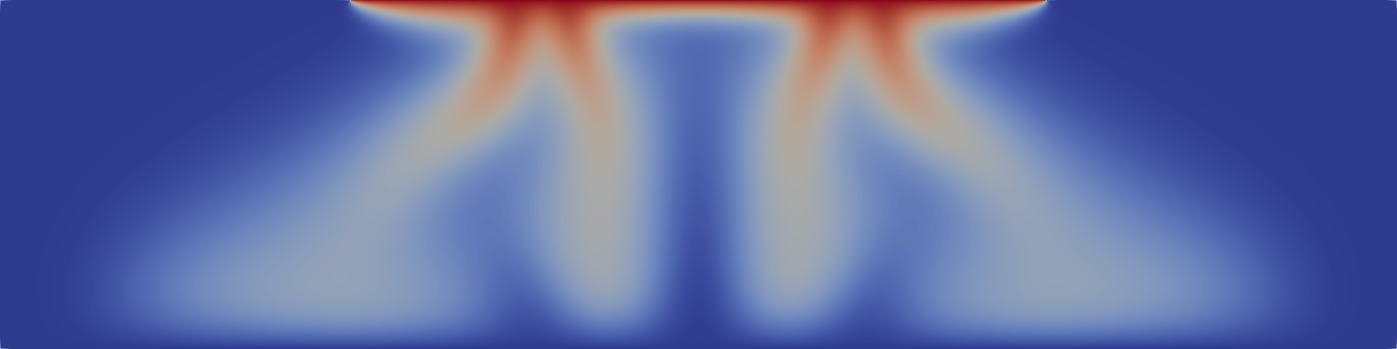}}
\caption{Solution $\sol$ at time $T=7$ years computed on the grids (a) of $16{,}384$ quadrilaterals with time step $0.007$ year and (b) of $262{,}144$ quadrilaterals with time step $0.00175$ year; $\poro=0.1$}
\label{fig:Ex1-grid-conv}
\end{figure}

The considered Elder's problem is known to be sensitive to the spatial and temporal grid resolution. Choice of a too coarse grid may produce an essentially different numerical solution. For that,
before we start to compare results of \gPC and \QMC, we should check if the chosen spatial mesh and the time step $\Delta t$ are sufficiently fine. Referring to the porosity in Eq.~\ref{eq:poro3},
we take $\poro(\bx,\thetab)= 0.1$, and compute the solution $\sol$ at time $T=7$ years on two different meshes: the first with $[N_0$ quadrilaterals, $\Delta t]$ and the second with $[16N_0$ quadrilaterals, $\frac{1}{4}\Delta t]$. The second grid is obtained from the first by two regular refinements. The results are presented in Fig.~\ref{fig:Ex1-grid-conv}. Both pictures looks very similar, i.e. the solution $\sol$, computed on coarse and fine meshes are almost equal. Although this fact is verified only for one realization, we expect the same result for the other sufficiently smooth porosities in this domain. Thus, the chosen spacial $N_0=16{,}384$ and temporal $\Delta t=0.007$ year resolutions are assumed to be sufficient for the following experiments with uncertain porosity.
\subsubsection{Comparison of \gPC and \QMC}
Now, we compare solutions computed by quasi Monte Carlo and \gPC methods. The number of the \QMC simulations is 600. The \gPC polynomial order is $p=5$. The 69 Clenshaw-Curtis quadrature points were used to compute \gPC coefficients.
Every simulation have been computed on a structured grid of $16{,}384$ quadrilaterals, there were $33{,}410$ degrees of freedom per realization. The time step was $\Delta t=0.007$ year. 

Every realization was computed in parallel on a separate node. Thus, to compute 600 \QMC simulations we used 600 parallel nodes. Each node consists of 32 parallel cores. The computing time of one realization is 1-2 minutes. The total computing time of all realizations is also about 1-2 minutes (due to parallelism).

\begin{figure}[htbp!]
\centering
\subfloat[]{\includegraphics[width=0.49\textwidth]{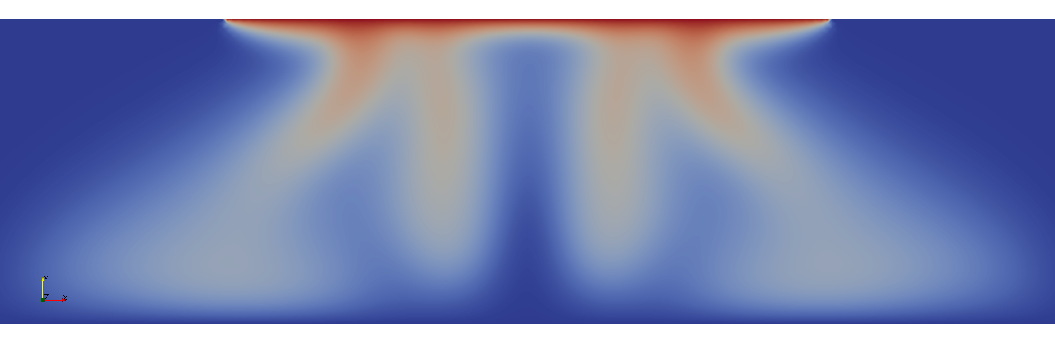}}\,
\subfloat[]{\includegraphics[width=0.49\textwidth]{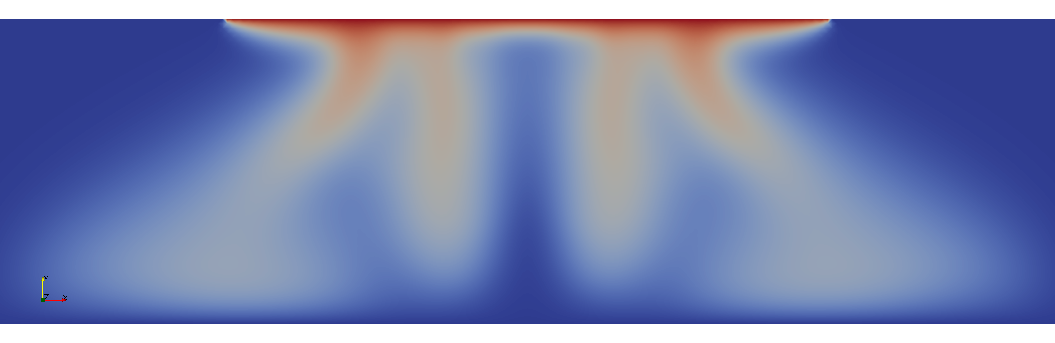}}\\
\subfloat[]{\includegraphics[width=0.49\textwidth]{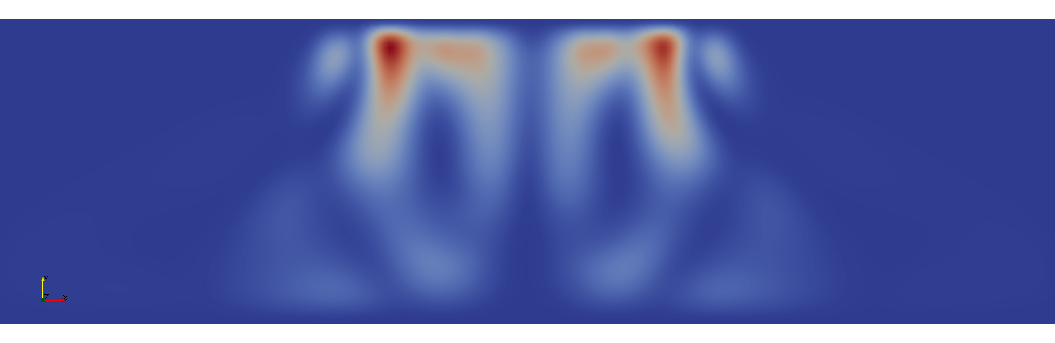}}\,
\subfloat[]{\includegraphics[width=0.49\textwidth]{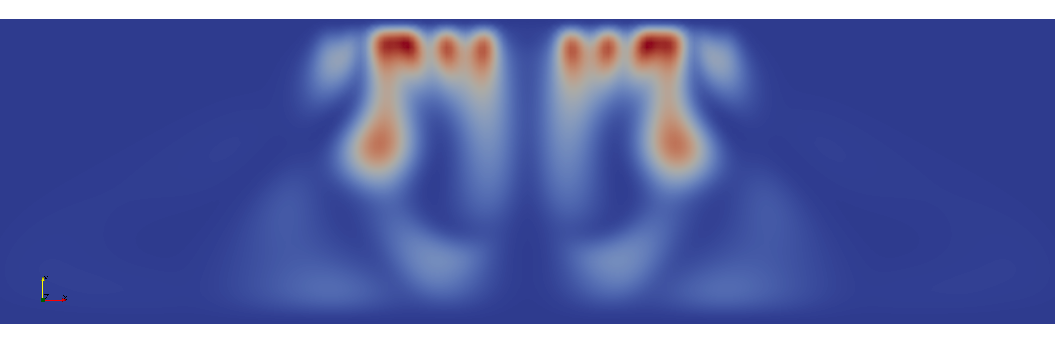}}
\caption{(a)-(b) The mean values of $\sol$ computed at time $T=7$ years by 600 \QMC and \gPC ($m=3$, $p=5$) respectively. Both methods give $\overline{\sol(\bx)}\in (0,1)$. (c)-(d) The variances of $\sol$ computed by the same methods, $\var{\sol}_{\mbox{\QMC}}\in (0,0.0498)$, and  $\var{\sol}_{\mbox{\gPC}}\in (0,0.0439)$; $p=5$, $m=3$.}
\label{fig:meanvar3rv}
\end{figure}

In Fig.~\ref{fig:meanvar3rv}(a)--(b) the mean values computed by \QMC and \gPC, respectively, at time $T=7$ years (1000 time steps for every realization) are shown.
Both mean values $\overline{\sol}_{\QMC}(\bx),\overline{\sol}_{gPC}(\bx)\in [0,1]$ are very similar. The variances are also computed by \QMC and \gPC and shown in Fig.~\ref{fig:meanvar3rv}(c)-(d). Both variances are very similar too:
$\var{\sol}_{\mbox{\QMC}}\in (0,0.0498)$, $\var{\sol}_{\mbox{gPC}}\in (0,0.0439)$. However, there are differences in the details, which can be explained by the approximation and integration errors.
\newpage
\subsection{Experiment 2. Elder's problem with 5 RVs and three layers}
The further example has more realistic settings. 
We consider the problem from Section \ref{litv:subsec:setup} (cf.\ Fig.\ \ref{fig:Elder2}) with three hydrogeological layers. The porosity field depends on 5 random variables:
\begin{equation}
\label{eq:poro5rvs}
\poro(\bx,\omega)=
\Bigg\{
 \begin{array}{cc}
0.08+0.01\sum_{i=1}^5 \theta_i \sin(ix\pi/600)\sin(iy\pi/150), & 120\leq y \leq 150.\\
0.06+0.01\sum_{i=1}^5 \theta_i \sin(ix\pi/600)\sin(iy\pi/150), & 50 \leq y <120 \\
0.09+0.01\sum_{i=1}^5 \theta_i \sin(ix\pi/600)\sin(iy\pi/150), & 0 \leq y < 50
\end{array}
\end{equation}
Here $y=0$ corresponds to the bottom and $y=150$ to the top. The porosity $\poro(\bx,\omega)$ has jumps between the layers. One of its realizations is shown in Fig.~\ref{fig:m5_3layers_sol}(e), where $\poro$ takes values in the interval $[0.0514, 0.09]$.

In this experiment, every realization has been computed on the regular grid of $1{,}048{,}576$ quadrilaterals with $2{,}102{,}274$ degrees of freedom. The time step is $\Delta t=0.005$ year. Every realization was computed in parallel on a separate node with 32 parallel cores. The computing time on each node is $~2.5$ hours.
%
\subsubsection{Dependence of the solution on the spatial and temporal grid}
\begin{figure}[htbp!]
\centering
\subfloat[]{\includegraphics[width=0.49\textwidth]{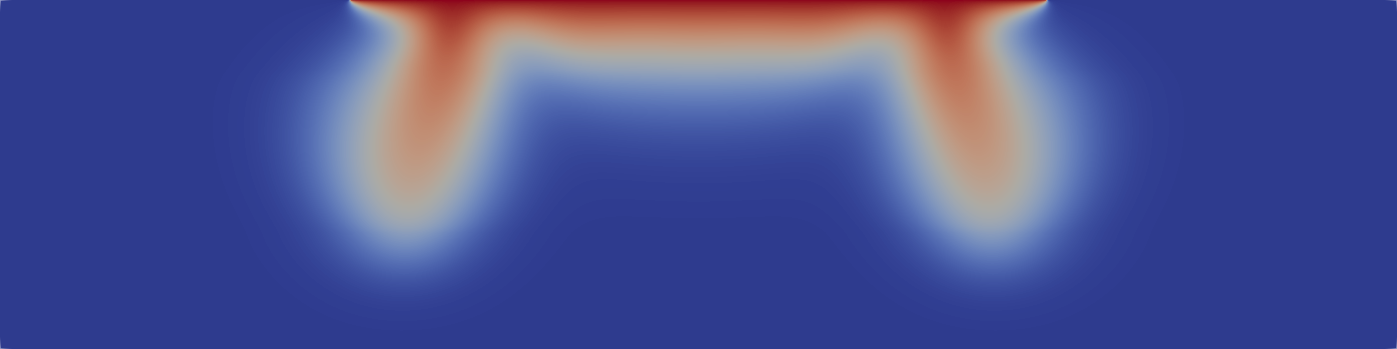}}\,
\subfloat[]{\includegraphics[width=0.49\textwidth]{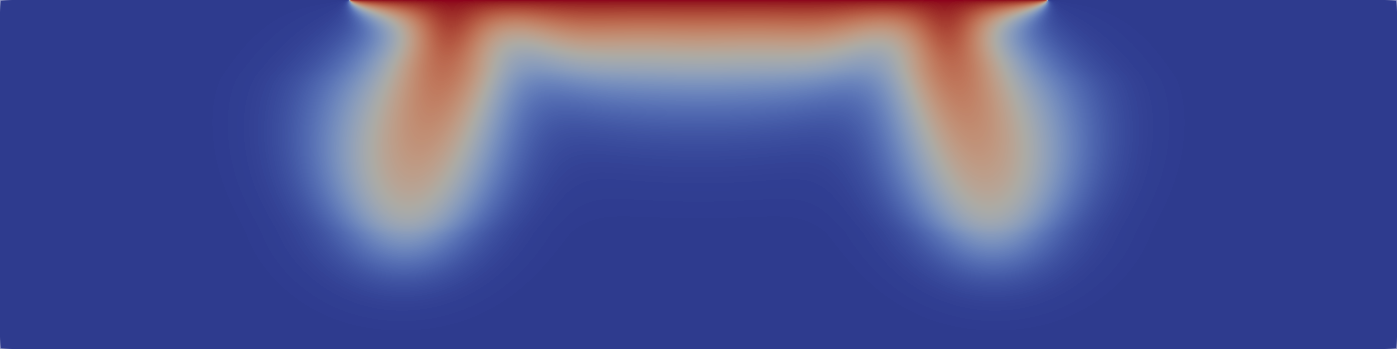}}
\caption{Solution $\sol$ at time $T=5.5$ years computed on the grids (a) of $1{,}048{,}576$ quadrilaterals with time step $0.005$ year and (b) of $4{,}194{,}304$ quadrilaterals with time step $0.0025$ year}
\label{fig:Ex2-grid-conv}
\end{figure}

Similar to the previous experiment, before we start to compare results of \gPC and \QMC, we check if the chosen spatial mesh and the time step $\Delta t$ are sufficiently fine.
We take $\poro(\bx,\thetab)= 0.1$, and compute the solution $\sol$ at time $T=5.5$ years on two different meshes: the first with $[N_0$ quadrilaterals, $\Delta t]$ and the second with $[4N_0$ quadrilaterals, $\frac{1}{2}\Delta t]$. The second grid is obtained from the first by one regular refinement. The results are presented in Fig.~\ref{fig:Ex2-grid-conv}. Both pictures looks very similar, i.e. the solution $\sol$, computed on coarse and fine meshes are almost equal. Thus, this experiment let us assume that $1{,}048{,}576$ quadrilaterals with the time step $0.005$ year are sufficient.
\subsubsection{Comparison of \gPC and \QMC}
\begin{figure}[ht]
   \centering
      \subfloat[]{\includegraphics[width=0.49\textwidth]{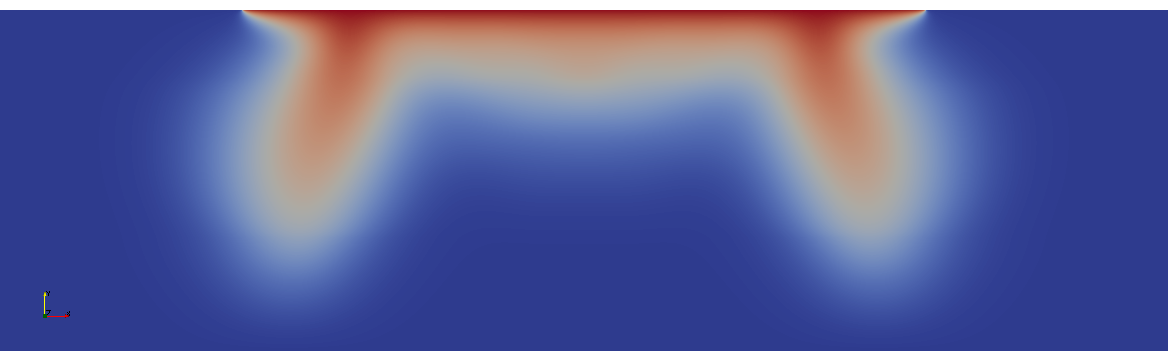}}\;
      \subfloat[]{\includegraphics[width=0.49\textwidth]{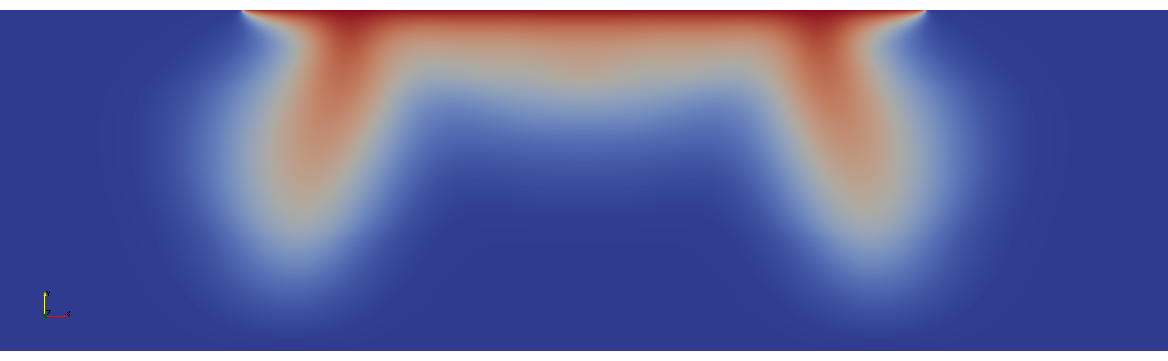}}\\
      \subfloat[]{\includegraphics[width=0.49\textwidth]{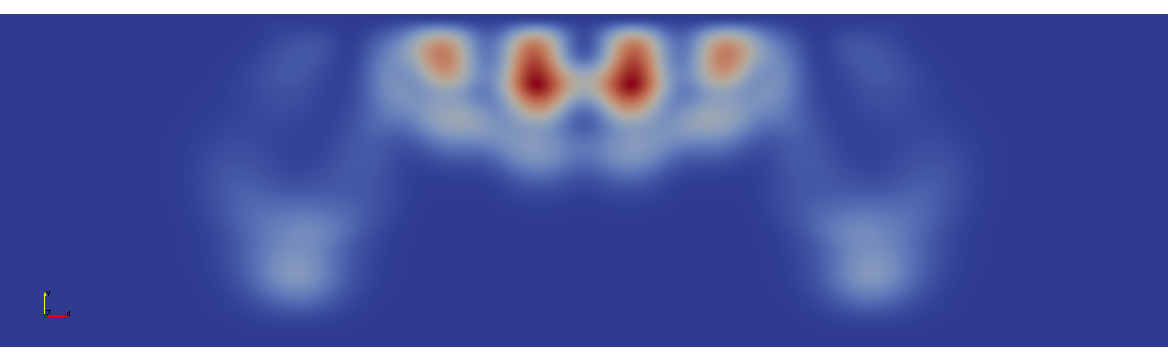}}\;
      \subfloat[]{\includegraphics[width=0.49\textwidth]{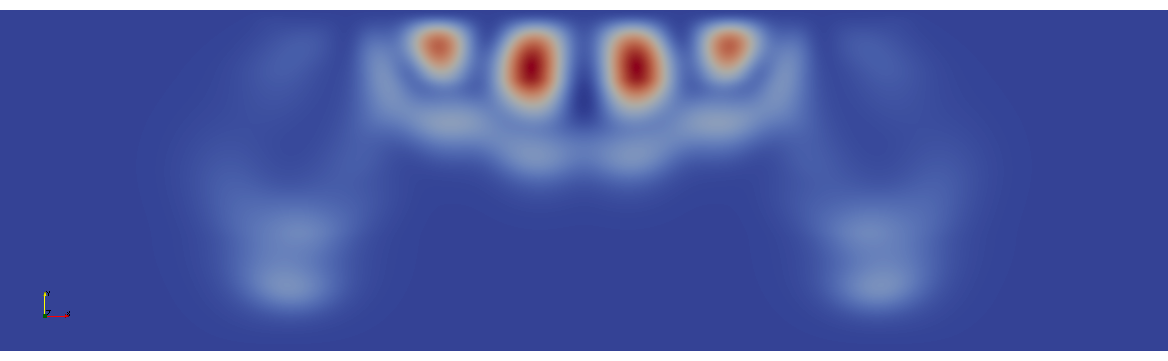}}\\
      \subfloat[]{\includegraphics[width=0.49\textwidth]{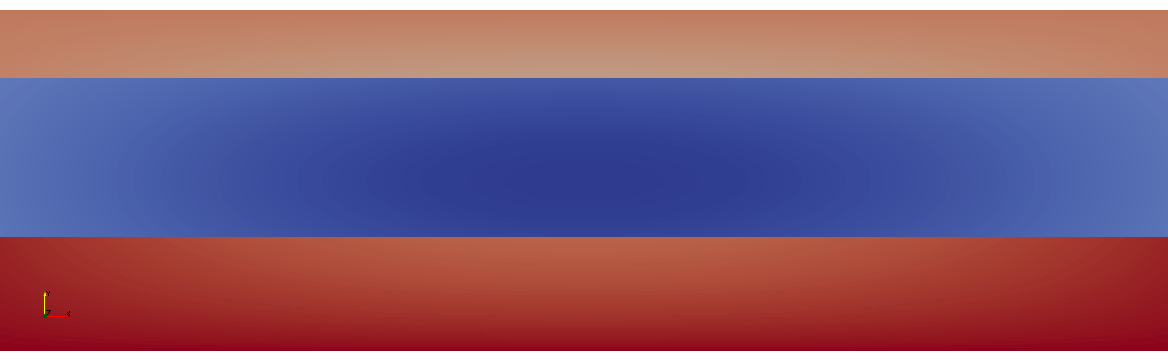}}\;
      \subfloat[]{\includegraphics[width=0.49\textwidth]{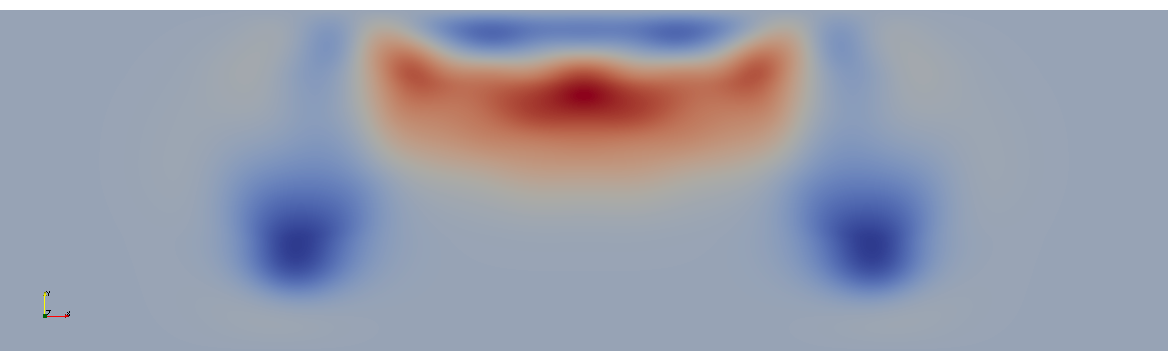}}
       \caption{a) $\overline{\sol}_{\gPC}(\bx,t)\in[0,1]$; b) $\overline{\sol}_{\QMC}(\bx,t)\in[0,1]$; c) $\var{\sol}_{\gPC}(\bx,t)\in[0,0.0466]$; d) $\var{\sol}_{\QMC}(\bx,t)\in[0,0.0556]$; e) three layers with different porosity, $\min_{\bx,\thetab}(\poro(\bx,\thetab))=0.0514$, $\max_{\bx,\thetab}(\poro(\bx,\thetab))=0.09$; f) difference between the deterministic solution (corresponds to $\btheta=\bzero$) and the mean value, $ \overline{\sol}(\bx,t)-\sol(\bx,t,\bzero)\in [-0.095,0.14]$ after $t=5.5$ years.}
\label{fig:m5_3layers_sol}
\end{figure}

We compare the mean values of $\sol$ computed by the \gPC (Fig.~\ref{fig:m5_3layers_sol}(a)) and by the \QMC method (Fig.~\ref{fig:m5_3layers_sol}(b)) at $T=5.5$ years (after 1100 time steps). We used 241 Clenshaw-Curtis quadrature points to compute \gPC coefficients.
Both mean values are very similar, and both values $\overline{\sol}_{\QMC}$ and $\overline{\sol}_{\gPC}$ variate in $[0,1]$. 
Figures~\ref{fig:m5_3layers_sol}(c)-(d) show the variances $\var{\sol}_{\gPC}(\bx,t)\in [0,0.0466]$ and $\var{\sol}_{\QMC}(\bx,t)\in [0,0.0556]$ respectively.
The difference $\sol(\bx,t,\bzero) - \overline{\sol}(\bx,t)$ between the deterministic solution (corresponding to $\btheta=\bzero$) and the mean value, computed for $t=5.5$ years, is presented in Fig.~\ref{fig:m5_3layers_sol}(f). This difference varies in the interval $[-0.095,0.14]$ and indicates the error which we find when we replace stochastic $\poro(\bx,\btheta)$ by $\poro(\bx,\btheta)|_{\btheta=\bzero}$. This large difference provides further indicator that it is necessary to estimate the propagation of uncertainties.  
\subsubsection{Computing probability density functions}
\label{sec:pdf}
Approximation of the probability density function of $\sol$ in a fixed point $(t^{*},\bx^{*})$ is computed by sampling the multivariate polynomial on the right-hand side in Eq.~\ref{eq:PCEpoint}. 
The random variable $\widehat{\mysol}(\thetab)$ can be sampled, e.g., $N_s=10^6$ times (no additional extensive simulations are required). The obtained samples can be used to evaluate the required statistics, and the probability density function.

\begin{figure}[htbp!]
  \centering
      \subfloat[]{
    \includegraphics[width=0.49\textwidth]{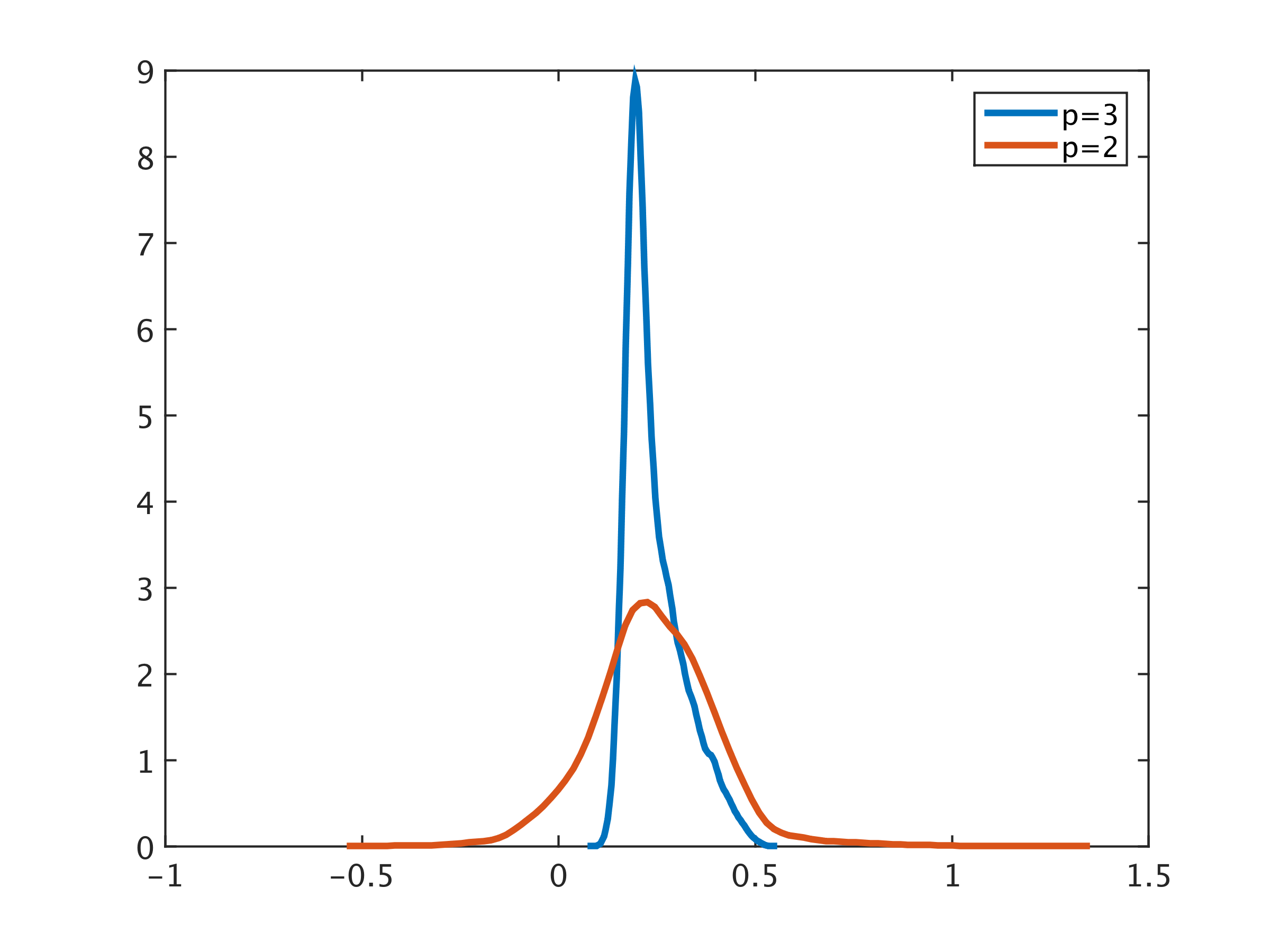}}\;
    \subfloat[]{\includegraphics[width=0.49\textwidth]{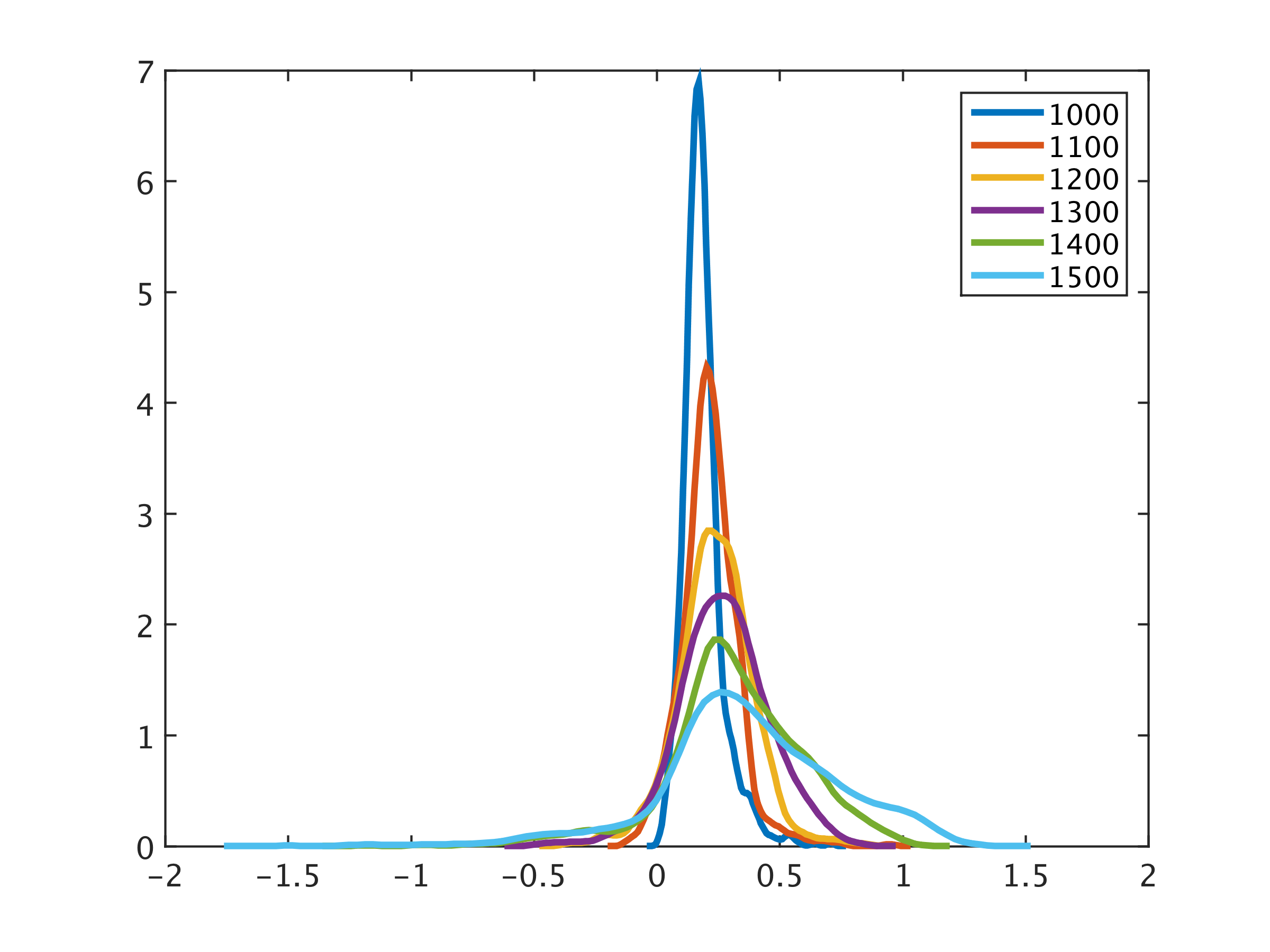}}\\
\caption{Comparison of pdfs at point $\bx=(300,100)$. a) pdfs sampled from \gPC of orders $p=2$ (blue curve) and $p=3$ (red curve). \gPC coefficients are computed with 241 quadrature points; b) pdfs computed after $\{1000,1100,\ldots, 1500\}$ time steps}
\label{fig:pdfs}
\end{figure}

Figure~\ref{fig:pdfs}(a) shows two pdfs computed with \gPC of orders $p=2$ and $3$, respectively, in point $\bx=(300,100)$ and after $1500$ time iterations (one iteration is 0.005 year). This figure shows that \gPC of order $p=3$ produces a pdf with a smaller variance. The mean values are equal. 
Figure~\ref{fig:pdfs}(b) shows six pdfs, computed with \gPC of orders $p=3$ in point $\bx=(300,100)$ and after $\{1000,1100,\ldots, 1500\}$ time steps or \textcolor{black}{$5.0-7.5$ years}. We can observe how the pdf function in that point is changing in time and that the variance is growing with time. 
The mean values are equal. 
\subsection{Computing quantiles}
\label{sec:quant}
We compute quantiles at three different points, e.g., $(150,50)$, $(300,75)$, $(500,100)\in \mathcal{D}$.
Table~\ref{tab:quantiles} presents the mean and the standard deviation for these three points after $1.375$, $2.75$ and $5.5$ years. We can see that the variance is increasing with increasing time steps (from bottom to top).

Below we calculate the quantiles for the cumulative probabilities $\{0.025,0.25,0.5,0.75,0.975\}$. For instance, the number $0.124$ in the first row and in the last column (Table~\ref{tab:quantiles}) means that $\mathbb{P}(\sol \leq 0.124)\geq 0.975$, i.e. the probability that the mass fraction $\sol\vert_{\bx=(150,50)}$ is smaller than $0.124$ is almost 1.
\begin{table}[ht!]
\caption{The mean, standard deviation (st.d.) and five quantiles at three points after $\{1.375, 2.75, 5.5\}$ years.}
\label{tab:quantiles}
\centering
\begin{tabular}{|c|c|c|c|c|c|c|c|c|} \hline
time, years   &point   &   mean & st.d. & \multicolumn{5}{c|}{five quantiles} \\ 
& $\bx=(x,y)$   &  & & $2.5\cdot 10^{-2}$& $2.5\cdot 10^{-1}$& $.50$& $.75$& $ .975$  \\ \hline
1.375&$(150,50)$  & $6.2\cdot 10^{-5} $& $2.8\cdot 10^{-5}$ & 0.021&   0.04&   0.058&   0.08& 0.124\\
2.75&$(150,50)$  & $3.2\cdot 10^{-1}$ & $2.2\cdot 10^{-2}$ &  0.273&   0.305&   0.322&   0.337& 0.356\\
5.5& $(150,50)$  & $3.8\cdot 10^{-1}$ & $2.2\cdot 10^{-2}$ & 0.339&   0.362&   0.379&   0.397& 0.420\\ \hline
%
1.375&$(300,75)$  & $1.3\cdot 10^{-4} $&$2.3\cdot 10^{-6}$ & 0.122&   0.124&   0.126&   0.128& 0.13\\
2.75&$(300,75)$  & $3.2\cdot 10^{-3}$ & $4.0 \cdot 10^{-4}$ & 0.0025&   0.0028&   0.0032&   0.0035& 0.0038\\
5.5& $(300,75)$  & $1.1\cdot 10^{-1}$ & $3.0 \cdot 10^{-3}$ & 0.099&   0.103&   0.106&   0.108& 0.111\\ \hline
1.375&$(500,100)$ & $2.2\cdot 10^{-5} $& $5.4\cdot 10^{-6}$ & 0.147&   0.18&   0.211&   0.253& 0.352\\ 
2.75&$(500,100)$ & $2.5\cdot 10^{-3}$ & $1.1 \cdot 10^{-3}$ & 0.001&   0.0016&   0.002&   0.003 &0.005\\ 
5.5& $(500,100)$ & $2.4\cdot 10^{-2}$ & $5.0 \cdot 10^{-3}$ & 0.016&   0.020&   0.024&   0.028& 0.035\\ \hline
\end{tabular}
\end{table}

\section{Conclusion}
\label{litv:sec:Conclusion}

We quantified the propagation of uncertainties in the density-driven groundwater flow problem with unknown porosity and permeability. Namely, we estimated the mean value, variance, exceedance probabilities and probability density function of the mass fraction. The considered problem is time-dependent and non-linear. The solution process requires fine resolution in time and space. Random perturbations in the porosity may change the stochastic solution, i.e., the number of fingers, their intensity, shape, propagation time and velocity. 

We applied the already known generalized polynomial chaos expansion (\gPC) method to the Elder’s problem in a 2D rectangular domain, where the porosity was modeled as a random field with 3-5 random variables. We demonstrate that there is a significant difference between the deterministic and stochastic problem setups and solutions. The results, obtained by the \gPC method were compared with the results, obtained by the \QMC method (Halton sequence). 

We demonstrate that by using the \gPC method, we can significantly reduce the computational cost of the classical (quasi) Monte Carlo method. Both methods for a moderate variance of the porosity give similar values for the mean and the variance. We have also demonstrated how to use the parallel multigrid solver \myug as a black-box solver in the stochastic framework. 

The practical novelty of this work in the implementation of the \gPC method on a distributed memory system, where all solutions and \gPC coefficients are distributed over multiple parallel nodes. The maximal number of parallel processing units was $600\times 32$, where 600 is the number of parallel nodes, and 32 is the number of computing cores on each node. We see great potential in using the \myug library as a black box solver for quantification of uncertainties. Our results are reproducible, and the code is available online on the GitHub repository\footnote{\url{https://github.com/UG4/ughub.wiki.git}}. 

The obtained results can be helpful to further understanding of the dispersion of seawater intrusions into coastal aquifers, radioactive waste disposal or contaminant plumes. 

In this work, we did not aim to determine the optimal parameters of the \gPC method. We also did not aim to compare different UQ methods. This point will be addressed in our future publications. As a next step, we plan to couple the \myug multi-grid solver with the Multi-Level Monte Carlo method. We intend to compute the majority of samples on low-cost coarse grids and just a few more on an expensive fine grid.

\begin{acknowledgement}
This work was supported by the King Abdullah University of Science and Technology (KAUST) and by the Alexander von Humboldt Foundation.
We used the resources of the Supercomputing Laboratory at KAUST, under the development project k1051. We would like to thank the KAUST core lab for the assistance with Shaheen II supercomputer.
We would also like to acknowledge developers of the \myug simulation framework from Frankfurt University.
\end{acknowledgement}
\begin{appendices}
\chapter{Legendre polynomials}
For unknown parameters which have uniform distribution, we employ multivariate Legendre polynomials as the basis for our \gPC.
Multivariate Legendre polynomials are defined on $[-1,1]^M$, where $M$ is the dimension of the stochastic space. The first few Legendre polynomials are 
\begin{equation}
\label{eq:LegendreMonom}
\pol_0(x)=1,\quad \pol_1(x)=x,\quad \pol_2(x)=\frac{3}{2}x^2-\frac{1}{2},
\pol_3(x)=\frac{5}{2}x^2-\frac{3}{2}x,\quad 
\pol_4(x)=\frac{35}{8}x^4-\frac{30}{8}x^2+\frac{3}{8}.
\end{equation}
The Legendre polynomials are orthogonal with respect to the $L_2$ norm on the interval $[-1,1]$.
The recursive formula is
\begin{equation}
\pol_0(x)=1, \quad \pol_1(x)=x,\quad (n+1)\pol_{n+1}(x)=(2n+1)x\pol_n(x) - n\pol_{n-1}(x).
\end{equation}
The scalar product is 
\begin{equation}
\int_{-1}^1 \pol_n(x)\pol_m(x)\mu(x)dx=\frac{1}{2n+1}\delta_{nm},
\\;\text{or}\;
or
\int_{-1}^1 \pol_n(x)\pol_m(x)dx=\frac{2}{2n+1}\delta_{nm},
\end{equation}
where the density function $\mu(x)$ is a constant (is equal 1/2 on $[-1,1]$), and where $\delta_{nm}$ denotes the Kronecker delta. 
The $L^2[-1,1]$ norm $\Vert \pol_n\Vert =\sqrt{\frac{2}{2n+1}}$ and the weighted $L^{2}_{\mu}[-1,1]$ norm is $\Vert \pol_n \Vert =\sqrt{\frac{1}{2n+1}}$. The \gPC coefficients are
\begin{equation}
c_{\alpha}(x):=\frac{(c(\thetab),\Pol_{\alpha}(\thetab))_{L^2_{\mu}}}{\Vert \Pol_{\alpha} \Vert^2_{L^2_{\mu}}}=
\frac{\int_{-1}^1 c\Pol_{\alpha} \mu(\thetab) d\thetab}{\int_{-1}^1 \Pol_{\alpha}\Pol_{\alpha} \mu(\thetab) d\thetab}=
\frac{\mu\int_{-1}^1 c\Pol_{\alpha} d\thetab}{\mu\int_{-1}^1 \Pol_{\alpha}\Pol_{\alpha} d\thetab}=\frac{0.5\cdot \sum_{i=1}^{n_q} c(\thetab_i)\Pol_{\alpha}(\thetab_i)w_i} {\frac{1}{2{\alpha}+1}}
\end{equation}
And for $d>1$ we have 
\begin{equation}
c_{\alpha}(x):=\frac{(c(\thetab),\Pol_{\alpha}(\thetab))_{L^2_{\mu}}}{\Vert \Pol_{\alpha} \Vert^2_{L^2_{\mu}}}=
\frac{\int_{[-1,1]^d} c\Pol_{\alpha} \mu(\thetab) d\thetab}{\int_{[-1,1]^d} \Pol_{\alpha}\Pol_{\alpha} \mu(\thetab) d\thetab}=\frac{0.5^d\cdot \sum_{i=1}^{n_q} c(\thetab_i)\Pol_{\alpha}(\thetab_i)w_i} {\prod_{i=1}^d\frac{1}{2\alpha_i+1}}
\end{equation}
\chapter{Difficulties in computing statistics}
We consider the following uncertain porosity (input parameter)
\begin{equation}
\label{eq:diff_finders}
\poro(\bx, \thetab)=0.1+0.01\cdot(\theta_1\cos(x/1200)+\theta_2\sin(y/300)+\theta_3\sin(x/2400)),\quad x\in [0,600],\;y\in[0,150].
\end{equation}

\begin{figure}[htbp!]
    \centering
    \subfloat[]{\includegraphics[width=0.49\textwidth]{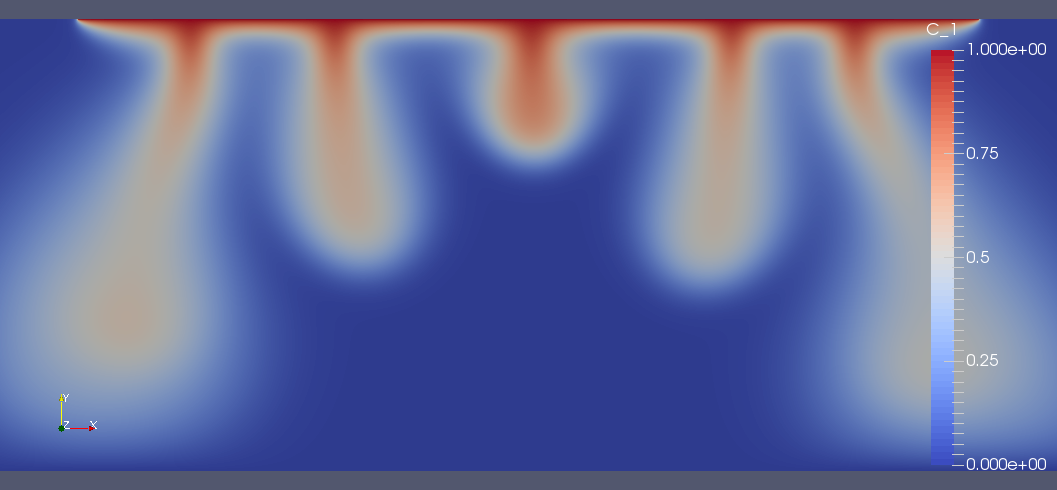}}\,
   \subfloat[]{    \includegraphics[width=0.49\textwidth]{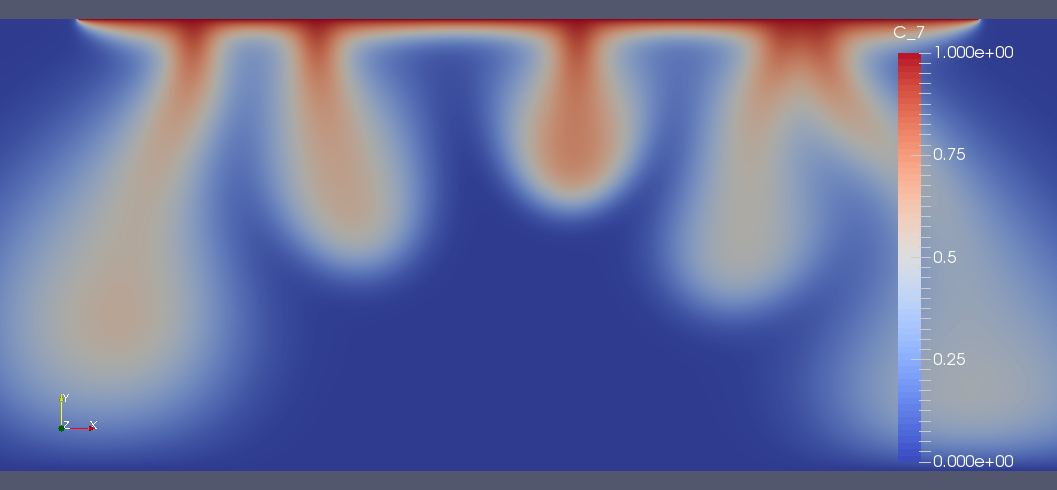}}\\
   \subfloat[]{    \includegraphics[width=0.49\textwidth]{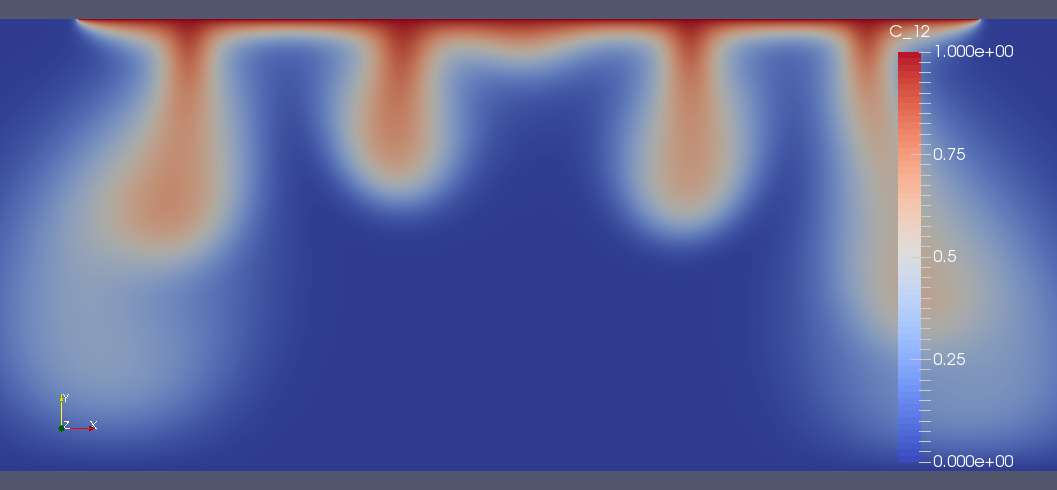}}\,
   \subfloat[]{    \includegraphics[width=0.49\textwidth]{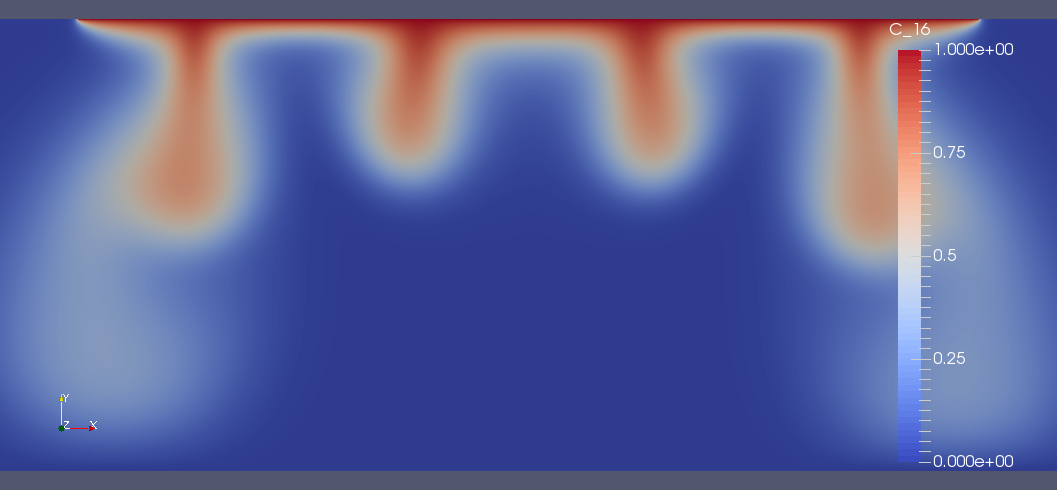}}   
\caption{Four different realizations of the mass fraction: with (a) 5, (b) 4.5, (c) 4, (d) 4 fingers. All $\sol\in[0,1]$.}
\label{fig:variances1-12}
\end{figure}

Figure~\ref{fig:variances1-12}
shows four different realizations of the solution (the mass fraction). These realisations correspond to 1st, 7th, 12th and 16th quadrature points. The number of fingers, their shapes, location and propagation are different. This fact makes it difficult to compute and to interpret the mean value, the variance and other statistics. 

\begin{figure}[htbp!]
    \centering
    \subfloat[]{\includegraphics[width=0.49\textwidth]{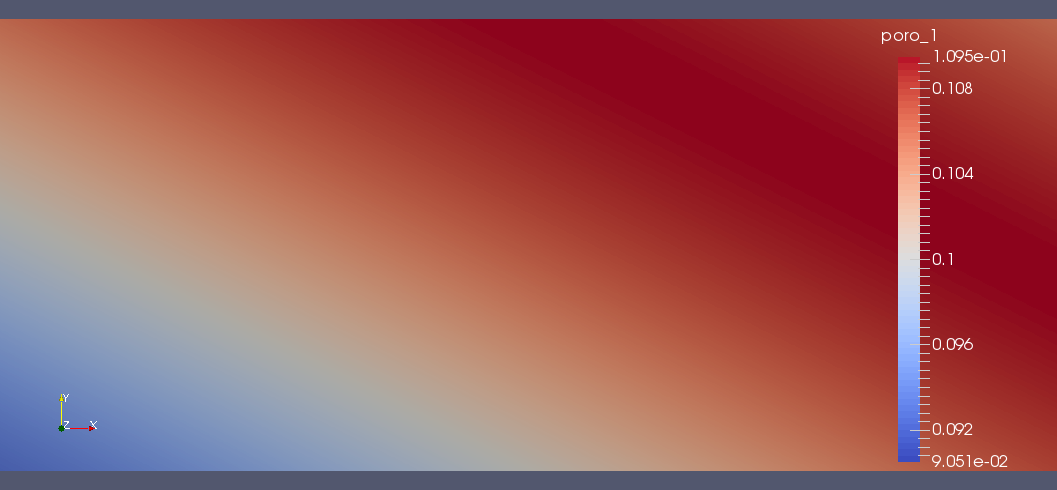}}\,
   \subfloat[]{    \includegraphics[width=0.49\textwidth]{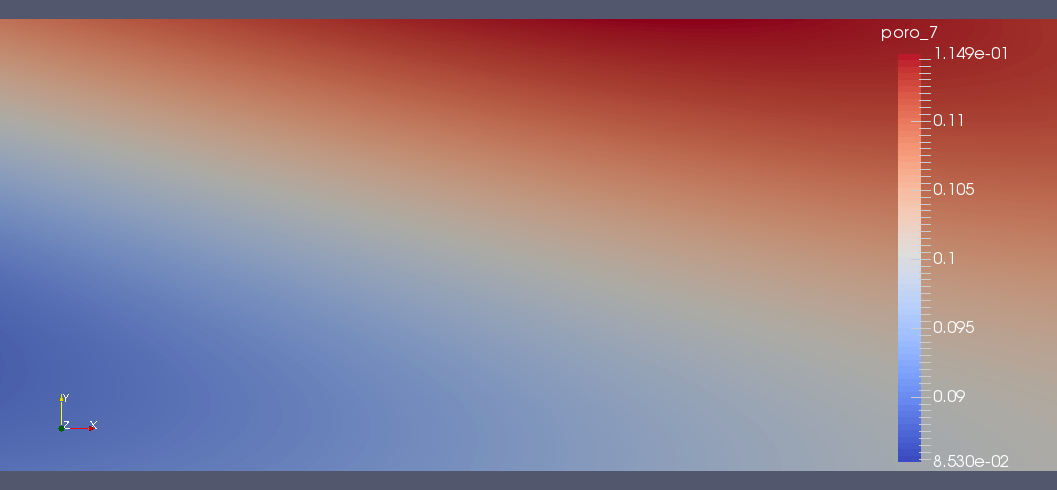}}\\
  \subfloat[]{ \includegraphics[width=0.49\textwidth]{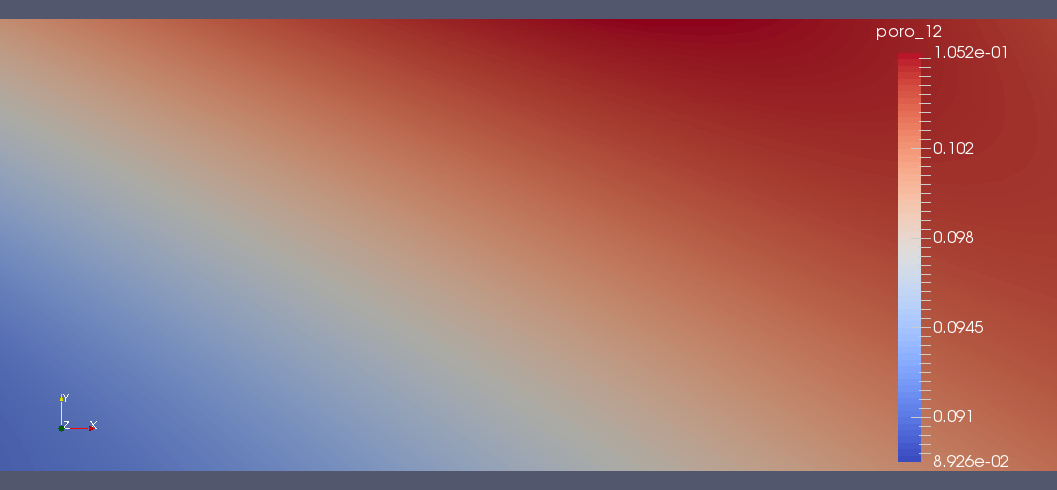}}\,
  \subfloat[]{ \includegraphics[width=0.49\textwidth]{poro7_m3_p7_small.png}}
\caption{Four realisations of the porosity field in Eq.~\ref{eq:diff_finders}.} 
\label{fig:poro_1_7}
\end{figure}

The porosity fields, corresponding to Fig.~\ref{fig:variances1-12}(a)-(b) are shown in Fig.~\ref{fig:poro_1_7}(a)-(b).

\end{appendices}


\end{document}